\documentclass[11pt]{article}
\usepackage{amsmath, amssymb}
\usepackage{amsthm}
\usepackage{amscd}
\usepackage{graphicx} 

\newcommand{\R}{\mathbb{R}}
\newcommand{\Z}{\mathbb{Z}}
\newcommand{\C}{\mathbb{C}}
\newcommand{\N}{\mathbb{N}}
\newcommand{\abs}[1]{\lvert#1\rvert}
\newcommand{\del}{\partial}

\newcommand{\restr}[1]{\bigr\rvert_{#1}}
\newcommand{\Lie}[1]{\mathfrak{#1}}	

\newcommand{\tensor}{\otimes}
\newcommand{\cross}{\times}

\newcommand{\into}{\hookrightarrow}
\newcommand{\define}[1]{\textsl{#1}}

\newcommand{\gap}{\vspace{2ex}}

\newcommand{\qn}{quantization}
\newcommand{\pol}{polarization}
\newcommand{\poln}{polarization}
\newcommand{\Q}{\mathcal{Q}}

\newcommand{\BS}{Bohr-Sommerfeld}
\newcommand{\lf}{\ell} 
\newcommand{\LL}{\mathbb{L}}	
\newcommand{\F}{\mathcal{F}}	
\newcommand{\sh}{\mathcal{S}}	
\newcommand{\J}{\mathcal{J}}	
\newcommand{\A}{\mathcal{A}}	
\newcommand{\B}{\mathcal{B}}  
\newcommand{\E}{\mathcal{E}}

\newcommand{\w}{\omega}	  
\newcommand{\T}{\mathcal{T}}	  
\newcommand{\ddtheta}{\frac{\del}{\del\theta}}       
\newcommand{\ddphi}{\frac{\del}{\del\phi}}
\newcommand{\ddth}[1]{\frac{\del}{\del\theta_{#1}}}  
\newcommand{\ddph}[1]{\frac{\del}{\del\phi_{#1}}}    
\newcommand{\bt}{\mathbf{t}}                         %
\newcommand{\bth}{\boldsymbol{\theta}}               %
\newcommand{\bs}{\mathbf{s}}                         %
\newcommand{\bph}{\boldsymbol{\phi}}                 %
\newcommand{\triv}{trivialization}  
\newcommand{\nbhd}{neighbourhood}
\newcommand{\pb}{^*\!}  
\newcommand{\id}{\text{id}}		
\newcommand{\msp}{(\R\cross S^1)^{n-k}\cross\C^k}	
\newcommand{\ip}[2]{\langle#1,#2\rangle}	
\newcommand{\inv}{^{-1}} 

\theoremstyle{definition}
\newtheorem{thm}{Theorem}[section]
\newtheorem{prop}[thm]{Proposition}

\newtheorem{lemma}[thm]{Lemma}
\newtheorem{cor}[thm]{Corollary}

\newtheorem{defnn}[thm]{Definition}
\newtheorem*{defn}{Definition}
\newtheorem*{notn}{Notation}	
\newtheorem*{rmk}{Remark}

\newtheorem*{thmm}{Main Theorem}
\newtheorem*{thmn}{Theorem}
\newtheorem*{lemn}{Lemma}

\theoremstyle{remark}
\newtheorem*{eg}{Example}

\newcommand{\epigraph}[2]{\begin{verse}\small{\it #1}\\ 
\hspace{5em} --- #2 \end{verse}}

\numberwithin{equation}{section}

 \newcommand     {\comment}[1]   {}
 \newcommand{\mute}[1] {}
 \newcommand     {\printname}[1] {}
\newcommand{\labell}[1] {\label{#1}\printname{#1}}

\begin{document}

\title{Locally toric manifolds and singular Bohr-Sommerfeld leaves}

\date{25 September, 2007}

\author{Mark Hamilton\thanks{Supported by a PIMS Postdoctoral Fellowship} \\
\texttt{umbra@math.utoronto.ca}}

\maketitle

\begin{abstract}

When geometric quantization is applied to a manifold using a real
polarization which is ``nice enough'', a result of \'Sniatycki says that the
quantization can be found by counting certain objects, called
Bohr-Sommerfeld leaves.  Subsequently, several authors have taken this as
motivation for counting Bohr-Sommerfeld leaves when studying the
quantization of manifolds which are less ``nice''.

In this paper, we examine the quantization of compact symplectic
manifolds that can locally be modelled by a toric manifold, using a real
polarization modelled on fibres of the moment map.  We compute the results
directly, and obtain a theorem similar to \'Sniatycki's, which gives the
quantization in terms of counting Bohr-Sommerfeld leaves.  However, the
count does not include the Bohr-Sommerfeld leaves which are singular.  Thus
the quantization obtained is different from the quantization obtained using
a K\"ahler polarization.  
\end{abstract}

\section{Introduction}\labell{s:intro}

Broadly speaking, \emph{quantization} is a procedure which associates 
to a symplectic manifold $M$ a Hilbert space $\Q(M)$.  
There are numerous methods of \qn; in this paper we consider
\emph{ geometric \qn.}
The ingredients for geometric \qn\ are: a symplectic manifold $(M, \omega)$, 
a complex line bundle $\LL$ over $M$, and a connection $\nabla$ on $\LL$
whose curvature is $\omega$.  The Hilbert spaces are constructed from
sections of $\LL$, or, more generally, higher cohomology groups associated 
with $\LL$.

One additional piece of structure is required, called a \emph{polarization.}
This is a method for restricting which sections of $\LL$ are considered; 
it is necessary because the space of all sections is usually ``too big.''
One example is a K\"ahler \poln, which is given by a complex structure
on $M$; the \qn\ is then determined in terms of holomorphic sections of $\LL$.
Another example is a \emph{real \poln,} which is defined more fully below.
It is then a natural question to ask, if we have different \poln s 
on a manifold $M$, do we obtain the same \qn s from them?  
This question can be called ``independence of \poln.''

A real \poln\ is given by a foliation of 
$M$ into Lagrangian submanifolds.  The sections of interest are those which
are \emph{leafwise flat:} 
covariant constant (with respect to $\nabla$) in the directions tangent
to the leaves of the foliation.  If $\J$ is the sheaf of leafwise flat
sections, then the \qn\ is constructed from the cohomology groups 
$ H^k(M;\J)$.

If the leaf space $B^n$ is a Hausdorff manifold and
the map $\pi \colon M^{2n} \to B^n$ is a fibration with compact fibres,
 a theorem of \'Sniatycki~\cite{Sn}
says that the above cohomology groups are zero except in dimension $n$; 
furthermore, $H^n$ can be expressed in terms of Bohr-Sommerfeld leaves.
A \emph{Bohr-Sommerfeld} leaf is one on which is defined a global
section which is flat along the leaf.
The set of \BS\ leaves is discrete, and \'Sniatycki's
theorem says that the dimension of $H^n$ is equal to the number of 
Bohr-Sommerfeld leaves.  (A similar result holds if the fibres are not 
compact, except the nonzero group occurs in a different dimension, equal
to the rank of the fundamental group of the fibres.)

Quantization using real polarizations, and the relation to the 
Bohr-Sommerfeld leaves, has also been studied by Andersen~\cite{andersen}.
He uses a different approach than we do in this paper, 
looking at the index theory of a ``leafwise complex.''
He generalizes \'Sniatycki's theorem to regular polarizations 
which are not necessarily real or K\"ahler, but something in between.

In many examples of interest, however, what seems like 
a real \poln\ is not quite a 
fibration, but can be viewed 
as a real \poln\ with singularities.  
Several authors, motivated by 
\'Sniatycki's theorem, have defined the \qn\ in such cases 
to be that obtained by simply counting 
Bohr-Sommerfeld leaves.  The characterization of \BS\ leaves often 
includes fibres which are singular; common practice has been to include
the singular fibres
 in the count, since  in many cases 
this gives a result which agrees with the
\qn\ obtained using a K\"ahler \poln .
(Two examples are Guillemin-Sternberg studying the Gelfand-Cetlin 
system in~\cite{GS},
and Jeffrey-Weitsman studying the moduli space of flat SU(2) connections
on a 2-manifold in \cite{JW}.)

One example of a ``singular real \poln '' is 
the level sets of the moment map on a toric manifold.
In this paper, we calculate directly 
the sheaf cohomology of a toric manifold, or 
more generally a manifold equipped with a Lagrangian fibration with elliptic
singularities, with respect to this singular polarization.
The result we obtain is similar to \'Sniatycki's 
theorem: all cohomology groups are zero except in one dimension, and 
the nonzero group can be computed by counting \BS\ fibres.  However, 
the singular \BS\ fibres are \emph{not} included in this count.

\begin{thmm}[\ref{mainthm}]
Let $M$ be a compact symplectic $2n$-manifold equipped with a 
locally toric singular Lagrangian fibration, with prequantization line bundle
$(\LL, \nabla)$.  Let $\J$ be the sheaf of leafwise flat sections 
of $\LL$.  Then the cohomology groups $H^k(M;\J)$ are zero for 
all $k \neq n$, and 
\[ H^n (M;\J) \cong \bigoplus_{b \in BS} \C \]
where the sum is taken over all nonsingular \BS\ fibres.
\end{thmm}

This is a surprising result, 
and is contrary to expectations.
It implies that this \qn\ is different from that computed using a 
K\"ahler \poln.  
For a toric manifold foliated by fibres of the moment map, 
the \BS\ leaves correspond 
to the integer lattice points in the moment polytope.
The \qn\ coming from a K\"ahler \poln\ has dimension equal to the number 
of lattice points in the polytope, including the ones on the boundary.
The \qn\ computed using the methods of this paper has dimension equal
to the number of lattice points in the \emph{interior} of the polytope, 
i.e.\ excluding the ones on the boundary.

The key point in the calculations where this difference appears 
is Theorem~\ref{disczero}, where 
we calculate the sheaf cohomology of a small neighbourhood 
of a singularity and find that it is zero. 
Thus, the singular leaf does not make a contribution to the 
cohomology, even though it is in some sense a Bohr-Sommerfeld leaf.  
This calculation relies on the observation in 
Proposition~\ref{prop:zero} that 
there are no nonzero leafwise flat sections defined on 
a neighbourhood of the origin.
The underlying reason for this lack of contribution, however, is unclear.


\subsection{Methods}

The method of proof we use is to compare the manifold to a standard 
model space.  We prove results for the model space by hand, and then
apply them locally to the manifold.  The types of spaces we consider 
are compact manifolds which, roughly, locally look like toric manifolds 
foliated by leaves of a moment map.  (See Definition \ref{spacedef} 
for the precise, if technical, statement.)  
These include not only toric manifolds, but also integrable systems 
with elliptic singularities.
We consider two types of model spaces:
the cylinder $\R \cross S^1$, and 
the complex plane $\C$.  

We equip each model space with a standard 
prequantization line bundle, 
and calculate the sheaf cohomology of each by hand 
(in sections \ref{s:cylinder} and \ref{C-section}).
We show how our methods apply to a simple example, that of $S^2$ 
acted on by the circle, in section~\ref{s:s2}.
After defining the model space for higher dimensions
in section~\ref{s:multdm}, 
in section \ref{s:leray}
we use a sheaf theoretic argument to patch together the lower-dimensional
results. 
The hypothesis that our manifold possesses a Lagrangian fibration with 
elliptic singularities guarantees that a neighbourhood of a \BS\ leaf
``looks like'' an open set in this model space, in a way that is 
compatible with the calculation of sheaf cohomology.  
In this way, we apply (in section \ref{s:piecing}) 
the results obtained for the model space to obtain the results for the 
general manifold.
Finally, we return to the specific case of toric manifolds in 
section~\ref{s:toric} to discuss the comparison between real and K\"ahler 
\poln s.

\gap
\textit{Acknowledgements:}
This paper is based on the author's Ph.D.\ thesis, 
carried out under the supervision
of Yael Karshon and Lisa Jeffrey.  
I am very grateful for all of their support over the course of my Ph.D.

I am grateful to Ruxandra Moraru for suggesting the methods of
Chapter~\ref{s:leray} and explaining some of the relevant sheaf
theory.  I also wish to thank Alejandro Uribe
for helpful comments; Eva Miranda, Paul Selick, and Pramath
Sastry for helpful suggestions and references; 
and Megumi Harada for much helpful advice.

\section{Background}\labell{s:bkgd}

\subsection{Connections}\labell{ss:connections}

Let $V$ be a vector bundle over a manifold $M$, 
$\Gamma(V)$ be the space of smooth\footnote{
In this paper, we always take sections, functions, and differential forms 
to be smooth.}
sections of $V$, 
and $\Omega^k(M)$ the space of (smooth) differential $k$-forms on $M$.
\begin{defn} \labell{def-conn}
Formally, a \define{connection} on a vector bundle $V$ is a map
$\nabla \colon \Gamma(V) \to \Omega^1(M) \tensor \Gamma(V)$
which satisfies the following properties:
\begin{enumerate}
\item $ \nabla (\sigma_1 + \sigma_2) = \nabla \sigma_1 + \nabla \sigma_2$
\item $ \nabla (f\sigma_1) = (df)\tensor \sigma_1 + f \nabla \sigma_1$
\end{enumerate}
for all sections $\sigma_1$ and $\sigma_2$ and functions $f$.
We typically write $\nabla_X \sigma$ for $\nabla \sigma$ applied to 
the vector (field) $X$.  
This is also called the \emph{covariant derivative} of $\sigma$ 
in the direction $X$.
\end{defn}

In order to express a connection in terms which are useful for 
calculations, we work locally. 
The following description, taken from \cite{Wood} (Appendix A.3),
only applies to (complex) line bundles, 
but that is all we require for this paper.

Let $\LL$ be a complex line bundle over some manifold $M$, 
with $s$ the unit section in some local \triv\ over $U\subset M$.
Fix a connection $\nabla$
on $\LL$.  Define the \define{potential one-form} $\Theta$ of the connection,
which is a one-form\footnote{If $\LL$ is a Hermitian line bundle,
and if $\nabla$ and $s$ respect the Hermitian structure, 
then $\Theta$ will in fact be real-valued.} 
on $U$, by 
\begin{equation}\labell{potl1form}
 \nabla_X s = -i\, \Theta(X) \, s. 
\end{equation}
The form $\Theta$ gives a complete description of the connection,
as follows: 
any other section $\sigma$ can be written as $\sigma = f s$ 
for some complex-valued function $f$, and so then using \eqref{def-conn},
we obtain
\begin{equation}\labell{eq:conn}
\nabla_X \sigma =df(X) s - f i\Theta(X) s.
\end{equation}
Usually, the trivializing section will be implicit, and we 
will blur the distinction between a section and the complex function
representing it.

{\bf Note:} This description of a connection
is only valid over an open set over which the bundle $\LL$ is trivializable.

\begin{defn}
The \define{curvature} of the connection 
is the 2-form $\Omega$ on $M$ defined by 
$\Omega = d\Theta$.
This is well-defined, independently of the choice of 
trivializing section (see Prop.~\ref{curvind}).
A connection is \define{flat} if its curvature is zero.
\end{defn}

The description of a connection in terms of its potential one-form 
depends on the choice of \triv.
In the following, we compute the effect that changing the trivialization
has on the potential one-form.

\begin{prop}\labell{chgtriv}
Let $\LL$ be a line bundle with connection $\nabla$ over $M$.
Suppose we have two local trivializations of $\LL$ over some open 
set, with unit sections $s$ and $s'$, respectively, which are 
related by $s=\psi s'$.  
(Here $\psi$ will be a nonvanishing $\C$-valued function.\footnote{
If the bundles are Hermitian, and the trivializations respect 
the Hermitian structure, $\psi$ takes values in $S^1$.}) 
If $\Theta$ and $\Theta'$ are the potential one-forms with respect 
to these trivializations, then
\begin{equation}
\Theta' = \Theta - i \tfrac{1}{\psi}\, d\psi.
\end{equation}
\end{prop}

\begin{proof}
Let $X$ be a vector.  Then, by eq.\ \eqref{potl1form},
\[ \nabla_X s = -i \Theta(X) \psi s'; \]
also, we have 
\[ \nabla_X \psi s' = d\psi(X) s' - i \Theta'(X) \psi s'. \]
Equating these and solving gives 
\[ -i \Theta(X) \psi s' = d\psi(X) s' - i\Theta'(X) \psi s'; \]
cancelling common factors (including $s'$)
and dividing by $\psi$ (which is never zero) gives 
\begin{equation}\labell{eq:trivchg}
\Theta' = \Theta - i \tfrac{1}{\psi}\, d\psi
\end{equation}
as one-forms.\footnote{
If we are in the Hermitian case, as in the previous footnotes, 
then $i \frac{1}{\psi} d\psi$ will be real-valued: \
since $\psi$ is $S^1$-valued,  
it can be written locally as $e^{if}$ for some real-valued function $f$, so
\[ d\psi = e^{if} i df. \]
Then
\[ i \tfrac{1}{\psi}\, d\psi = -df \]
is real-valued.\\
} 
\end{proof}

\begin{prop}\labell{curvind}
The curvature form $\Omega$ of a connection is independent of the 
trivialization used to compute the potential one-form.
\end{prop}

\begin{proof}
This is a straightforward calculation.  If we have two different 
trivializations, the potential one-forms are related by 
\[ \Theta' = \Theta - i\tfrac{1}{\psi} \, d\psi. \]
Thus 
\begin{equation}
\begin{split}
 \Omega' = d\Theta' &= d \bigl( \Theta - i\tfrac{1}{\psi}\, d\psi \bigr) \\
&= d\Theta - i \, d\bigl( \tfrac{1}{\psi}\bigr)  \wedge d\psi 
	+ i\tfrac{1}{\psi}\, d\, d\psi\\
&= d\Theta + i \tfrac{1}{\psi^2}\, d\psi \wedge d\psi\\
&= d\Theta.
\end{split}
\end{equation}
This also implies that $\Omega$ is defined on all of $M$.
\end{proof}

\subsubsection{Holonomy}\labell{s:hol}
Suppose we have a line bundle with connection $(\LL,\nabla)$.

Let $\gamma$ be a curve on $M$, with tangent vector $\dot\gamma$,
and suppose $\sigma$ is a section of $\LL$ defined at least on $\gamma$.
Let $\tilde\gamma$ be the lifting of $\gamma$ to $\LL$ via $\sigma$, 
i.e., $\tilde\gamma = \sigma \circ \gamma$.
Then the lifting $\tilde\gamma$ is \define{horizontal} if 
\begin{equation}
\nabla_{\dot\gamma}\sigma = 0
\end{equation}
for all points along the curve.

Given a curve $\gamma$ in $M$ and a point $x$ in the fibre 
over $\gamma(0)$, the curve $\tilde\gamma$ is uniquely determined
by the condition that it is a horizontal lift of $\gamma$ 
with $\tilde\gamma(0) = x$.
Taking $x$ to $\tilde\gamma(1)$ gives a map from $\LL_{\gamma(0)}$ to 
$\LL_{\gamma(1)}$, called \define{parallel transport} along $\gamma$.
If $\gamma$ is a loop, this gives an automorphism of 
$\LL_{\gamma(0)}$, called the \define{holonomy} around $\gamma$.

If $\LL$ is Hermitian and the connection respects the Hermitian structure, 
we can view holonomy as a map from \{loops on $M$\} to $S^1$.
This map is given by 
\begin{equation}\labell{holonomy}
\text{hol} \colon \gamma \mapsto \exp \left( i \int_\gamma \Theta \right).
\end{equation}
If the connection is flat, then homotopic loops have the same holonomy
(the proof is basically Stokes' theorem) 
and the holonomy can be viewed as a map
from $\pi_1(M)$ to the automorphisms of the fibre.

\subsection{Sheaves and cohomology}\labell{ss:sheaf-cohom}

We review the definition of sheaves and the construction of
\v Cech cohomology, as they will be used extensively in this paper.
The material in this section is standard, and so we do not prove
our assertions.  See, for example, \cite{GH}, pp.\ 38--40.

Roughly speaking, a sheaf is a collection of functions on open sets,
often satisfying some further conditions (holomorphic, etc.).
The precise definition is as follows:
\begin{defn}  Let $X$ be a topological space.  A \define{presheaf} 
(of abelian groups)  $\F$
on $X$ assigns to every open set $U$ of $X$ an abelian group $\F(U)$,
usually referred to as the set of \define{sections} of $\F$ over $U$.
It also assigns \define{restriction maps:}
to any $V \subset U$, the presheaf assigns a map $\F(U) \to \F(V)$,
``restriction to $V$'', such that if $W \subset V \subset U$ and 
$\sigma \in \F(U)$, then 
\[ \sigma\restr{W} = (\sigma\restr{V})\restr{W}, \]
and if $V=U$ then ``restriction'' is just the identity map.
\end{defn}

\begin{defn}\labell{defn-sheaf}
A presheaf $\sh$ is a \define{sheaf} if the following properties hold:
\begin{enumerate}
\item For any pair of open sets $U$, $V$, and sections $\sigma \in \sh(U)$
and $\tau \in \sh(V)$ which agree on the intersection $U\cap V$, 
 there exists a section $\rho \in \sh(U\cup V)$ which restricts to 
$\sigma$ on $U$ and $\tau$ on $V$.

\item If $\sigma$ and $\tau$ in $\sh(U\cup V)$
have equal restrictions to $U$ and $V$, then they are equal
on $U\cup V$.
\end{enumerate}
\end{defn}

\begin{defnn}\labell{d:stalk}
For a sheaf $\sh$ over $M$, the \define{stalk} 
of $\sh$ over $x \in M$ is 
\[ \sh_x =  \varinjlim \sh(U) \]
where the limit is taken over all open sets $U$ containing $x$.
\end{defnn}

The \define{\v Cech cohomology of $M$ with coefficients in the sheaf $\sh$}
(or just ``the sheaf cohomology of $M$'')\footnote{
The ``sheaf cohomology of $M$'' 
is actually defined in a more abstract way using resolutions 
of the sheaf $\sh$.
However, for manifolds, the result obtained is the same as the
\v Cech cohomology, and the latter is more convenient for our calculations.}
is defined as follows.

Fix an open cover $\A = \{ A_\alpha\}$ of $M$.

A \define{\v Cech $k$-cochain} assigns,
to each $(k+1)$-fold intersection
of elements from the cover $\A$, a section of $\sh$.  
We write $A_{\alpha_0 \cdots \alpha_k}$ for 
$A_{\alpha_0} \cap \cdots A_{\alpha_k}$,
where the $\alpha_j$ are distinct.  
Then a $k$-cochain is an assignment 
$f_{\alpha_0 \cdots \alpha_k} \in \sh(A_{\alpha_0 \cdots \alpha_k})$
for each $(k+1)$-fold intersection in the cover $\A$.
Denote the set of $k$-cochains by $C^k_\A(M;\sh)$, 
or just $C^k_\A$ when the manifold and sheaf are understood.

Next, we define a coboundary operator $\delta$ to make $C^\ast_\A$ into 
a cochain complex.  
For $f = \{f_{\alpha_0 \cdots \alpha_{k-1}} \}$ a $(k-1)$-cochain, 
$\delta f$ will be a $k$-cochain; thus, 
we need to give a section corresponding to each $(k+1)$-fold intersection.
This is done as follows:
\begin{equation}\labell{defdelta}
 (\delta f)_{\alpha_0 \cdots \alpha_k} = 
\sum_{j=0}^{k} (-1)^j f_{\alpha_0 \cdots \hat{\alpha}_j \cdots \alpha_k}
\restr{A_{\alpha_0 \cdots \cdots \alpha_k}}
\end{equation}
where the $\hat{}$ denotes that the index is omitted.
Thus, for instance, 
$(\delta f)_{123} = f_{23} - f_{13} + f_{12}$, all restricted to $A_{123}$.

A (straightforward, but tedious) calculation shows that 
$\delta\circ\delta=0$, 
and so $C^\ast_\A$ is a cochain complex.
The sheaf cohomology with respect to the cover $\A$ is the cohomology 
of this complex, 
\[ H^k_\A (M;\sh) = \frac{\ker \delta^{k}}{\text{im}\, \delta^{k-1}} \]
(where by $\delta^k$ we mean the map $\delta$ on $C^k_\A$).

Another cover $\B$ is a \define{refinement} of $\A$,
and we write $\B \leq \A$, 
if every element of $\B$ is a subset of some element of $\A$.  
From this we define a map 
$\phi \colon C^k_\A(U,\sh) \to C^k_\B(U,\sh)$, induced by the restriction 
maps in the sheaf: simply restrict each element of a cochain, 
defined on some intersection of sets in $\A$,
to the intersection of the corresponding sets in $\B$.
More formally, a refinement gives a map $\rho \colon \B \to \A$, 
where $B \subset \rho(B)$ for all $B\in \B$.
Then, if $\eta \in C^k_\A$ is a cochain, 
$\phi \eta$ is defined by 
\[ (\phi \eta)_{B_0 B_1 \cdots B_k} = 
(\eta)_{(\rho B_0)  (\rho B_1)  \cdots ( \rho B_k)}
\restr{B_0 B_1 \cdots B_k}. 
\]

Since this map is essentially just restriction, it commutes
with $\delta$, and so it induces a map on cohomology $H^*_\A \to H^*_\B$.
Two different choices of  maps $\rho$ for the same refinement induce
chain homotopic maps on cochains, and thus induce the same map on cohomology.
These maps turn the collection of $H^\ast_\A$ for all open 
covers of $M$ into a directed system. 

Finally, the honest-to-goodness \define{sheaf cohomology} 
of $M$ is defined as the limit of this directed system:
\[ H^*(M;\sh) = \varinjlim H^*_\A(M;\sh). \]
 
\begin{lemma}\labell{cofinal}
Suppose $\mathfrak{B}$ is a collection of covers $\B$ of $M$ 
such that any open cover $\A$ of $M$ has a refinement 
$\B \in \mathfrak{B}$.  Suppose furthermore that all 
$H^*_\B(M;\sh)$ are isomorphic, for each $\B \in \mathfrak{B}$.
Then the sheaf cohomology of $M$ is isomorphic to the cohomology
computed using one of the covers $\B$.
\end{lemma}

\begin{proof}
This follows easily from the construction of 
the direct limit.
The set $\mathfrak{B}$ is \define{cofinal} in the set of all covers if
every cover has a refinement from $\mathfrak{B}$.
\end{proof}

\subsubsection{Naturality of Direct Limits and Cohomology}\labell{ss:natural}

The following results will be used in Section \ref{ss:cyl-MV}.

\begin{lemma}[\cite{Gode}, p.\ 10]\labell{limnat}
The direct limit of exact sequences is exact.  
More precisely, suppose we have:
\begin{itemize}
\item  three directed systems
of modules $L_i'$, $L_i$, and $L_i''$ (indexed by the same directed set $I$), 
with direct limits $L'$, $L$, and $L''$ respectively, and

\item for each $i$, an exact sequence $L_i' \to L_i \to L_i''$.
\end{itemize}
Suppose further that, for each $i \geq j$, the diagram 
\[
\begin{CD}
	L_i'	@>>>	L_i	@>>>	L_i''\\
	@VVV		@VVV		@VVV\\
	L_j'	@>>>	L_j	@>>>	L_j''
\end{CD}
\]
where the vertical arrows denote the maps in the directed system,
is commutative.

Then the sequence $L' \to L \to L''$ is exact.
\end{lemma}

\begin{lemma}[\cite{MacLane}, Prop.\ 4.2]\labell{sesnat}
Given a transformation of short exact sequences of cochain complexes
\[ 
\begin{CD}
	0	@>>>	A^*	@>>>	B^*	@>>>	C^* @>>> 0\\
	@.		@VVV		@VVV		@VVV\\
	0	@>>>	D^*	@>>>	E^*	@>>>	F^* @>>> 0
\end{CD}
\]
this induces a map between the long exact sequences of cohomology
\[
\begin{CD}
\cdots	@>>>	H^*(A)	@>>>	H^*(B)	@>>> H^*(C)	@>>>	\cdots\\
@.		@VVV		@VVV		@VVV\\
\cdots	@>>>	H^*(D)	@>>>	H^*(E)	@>>> H^*(F)	@>>>	\cdots
\end{CD}
\]
\end{lemma}

\subsection{Toric manifolds}\labell{ss:toric}

We briefly review the definition of, and a few facts about, toric manifolds,
referring the reader to~\cite{ACL} for a more detailed introduction.

\begin{defn}
A \define{toric manifold} is a compact symplectic manifold $M$ of dimension 
$2n$ equipped with an effective Hamiltonian action of the torus 
$T^n = (S^1)^n$.
\end{defn}

Recall that if a symplectic manifold $M^{2n}$ has an 
effective Hamiltonian action of a torus $T^k$, then $k \leq n$.
Thus $M$ is toric if $T$ has maximal dimension.

The ``Hamiltonian'' in the definition means that 
a toric manfold is equipped with a $T$-invariant \define{moment map} 
$\mu \colon M \to \R^n \cong \Lie{t}^*$,
which satisfies the following condition:
For $\xi \in \Lie{t}$, let $\xi_M$ be the generating vector field for the 
action on $M$, and let $\mu^\xi\colon M \to \R$ be the 
``component of $\mu$ in the $\xi$ direction,'' that is, 
\[ \ip{\mu(x)}{\xi} = \mu^\xi(x)\]
for all $x \in M$,
where $\ip{\phantom{a}}{\phantom{a}}$ 
denotes the pairing between $\Lie{t}$ and $\Lie{t}^*$.
Then
\begin{equation}\labell{momap}
\imath_{\xi_M} \omega = d \mu^\xi, 
\end{equation}
that is, $\xi_M$ is the Hamiltonian vector field of $\mu^\xi$.
We do not use this definition directly,
but some facts about toric manifolds.

First, fibres of the moment map are orbits of the torus action.
Second, 
a theorem of Atiyah and Guillemin-Sternberg says that 
the image of the moment map is a convex polytope 
$\Delta$ in $\Lie{t}^* \cong \R^n$.
If $x\in\Delta$ lies on a face of codimension $k$,
then the stabilizer of a point in $\mu\inv(x)$
is a torus of dimension $k$, and $\mu\inv(x)$
is a torus of dimension $n-k$.
(This is Lemma 2.2 in~\cite{delzant}.)
For $x$ in the interior of $\Delta$,
 the fibres are $n$-dimensional tori; the moment map condition implies 
they are Lagrangian.

The following lemma, known as the Local Normal Form, gives a description 
of the neighbourhood of an orbit.

\begin{lemma}[Local Normal Form, Lemma 2.5 in \cite{delzant}]\labell{lnf}
Let $M$ be a symplectic toric manifold with moment polytope $\Delta$
and moment map $\mu$.
Let $F$ be a face of $\Delta$ of dimension $m$;
$V$ a convex, open, relatively compact subset of $F$; $x$ a point in $V$; 
and $D$ a ball around 0 in $\C^{2(n-m)}$.
Let $\omega_0$ be the following symplectic form on 
$(S^1)^m\cross V \cross D$:
\[ \omega_0 = \sum_{1\leq j\leq m} d\alpha_j \wedge da_j + 
\sum_{m+1 \leq j \leq n} dx_j \wedge dy_j \]
where $\alpha$ are coordinates on $S^1$, $a$ are coordinates on $V$, 
and $z=x+iy$ are coordinates on $D$.

Then there is a symplectic isomorphism from a neighbourhood of $\mu\inv(V)$
onto $(S^1)^m\cross V \cross D$ taking the action of $T$ to the 
action of $(S^1)^n$ defined by 
\begin{equation}
\begin{split}
(\theta_1, \ldots, \theta_n) \cdot 
(\alpha_1,\ldots,\alpha_m,a_1,\ldots,a_m, z_{m+1},\ldots,z_n) =\\
( \alpha_1+\theta_1 ,\ldots, \alpha_m+\theta_m, a_1,\ldots,a_m, 
e^{i\theta_{m+1}} z_{m+1}, \ldots, e^{i\theta_n} z_n) 
\end{split}
\end{equation}
with moment map 
\[ \mu = p + (a_1,\ldots,a_m,\abs{z_{m+1}}^2,\ldots, \abs{z_n}^2).\]
\end{lemma}

\subsection{Geometric \qn\ and polarizations}\labell{ss:qn-pol}

Let $(M,\w)$ be a compact symplectic manifold of dimension $2n$.

\subsubsection{Quantization}\labell{sss:qn}

The theory of geometric quantization was initiated by Kostant 
and Souriau in the 1970s, and remains an active area of research 
today, with applications to both physics and representation theory.
We do not attempt to give a comprehensive introduction here.  
For a nice, brief overview of the ideas behind it, 
see~\cite{ggk}, section 6.1, or~\cite{gsbook}, chapter 34.  
For a more thorough introduction, see~\cite{puta}.
Two classic references, albeit somewhat technical, are~\cite{sniat}
and~\cite{Wood}.

The basic idea of quantization is 
to associate to a symplectic manifold $(M,\w)$ a Hilbert space 
(or a vector space) $\Q(M)$.  
(The terminology ``quantization'' comes from physics, 
where we think of $M$ as 
a classical mechanical system, and $\Q(M)$ as the 
space of wave functions of the corresponding quantum system.)
Much of the motivation for geometric \qn\ in mathematics comes from 
representation theory.  

In geometric quantization, the quantum space is constructed from 
the sections of a complex line bundle or, more generally, 
from higher-dimensional cohomology groups associated with the line bundle.

\begin{defn}
A \define{prequantization line bundle} over $M$
is a Hermitian line bundle $\LL$ over $M$, 
with a connection $\nabla$ whose curvature 
is $\w$.
$M$ is \define{prequantizable} if it possesses a prequantization line bundle.
(This will be the case if{f} the symplectic form satisfies the integrality 
condition that $\frac{1}{2\pi} [\omega] \in H^2(m,\Z)$.  See, for example, 
\cite{Wood}, section 8.3.\footnote{
The exact form of the integrality condition depends on the conventions used,
and one may see instead
$\frac{1}{\hbar} [\omega]$, $\frac{1}{2\pi\hbar}[\omega]$, or just $[\omega]$
required to be integral.  We take $\hbar = 1$, and use coordinates on $S^1$ 
running from $0$ to $2\pi$, which gives the integrality 
condition stated here.})
\end{defn}

We would like $\Q(M)$ to be the space of sections of $\LL$.  
However, this space is generally ``too big.''
As noted in the Introduction, the solution is to use a 
``polarization'' to choose a subspace of the space of sections; 
the quantum space is then constructed from only ``polarized'' sections.

Our main interest in this paper is in \emph{real polarizations.}\footnote{
Thus in particular 
we do not give the general theory of polarizations, but refer 
the reader to~\cite{sniat}, pages 8--11, or~\cite{Wood} sections 4.5 and 5.4.
We will also have occasion to mention K\"ahler \poln s, which 
we define below.}
The usual definition of a
real polarization on $M$ is a sub-bundle $P \subset TM$ 
which is Lagrangian and integrable. 
In our case, we allow `singular polarizations', where the leaves
are not all of the same dimension. 
We define a locally toric singular Lagrangian fibration
to be a structure which locally looks like the 
(singular) fibration of a toric manifold by the moment map, 
using the local structure given by Lemma~\ref{lnf}.

\begin{defnn}\labell{spacedef}
A \define{locally toric singular Lagrangian fibration} 
on a symplectic $2n$-manifold $M$
is a map $\pi \colon M \to B$ 
to a topological space $B$ such that
for every point in $B$, there exist:
\begin{itemize}
\item a nonnegative integer $k$
\item a neighbourhood $U\subset B$ of the point
\item an open subset $\Omega \subset \R^{n-k} \cross \R_+^k$
\item a homeomorphism $\psi \colon U \cong \Omega$
\item a symplectomorphism 
$\widetilde\psi \colon \pi^{-1}(U) \cong \pi_0^{-1}(\Omega)$ 
\end{itemize}
such that the following diagram commutes:
\begin{equation}
\begin{CD}
\pi^{-1}(U)	@>\widetilde{\psi}>>	(S^1 \cross \R)^{n-k} \cross \C^k\\
\pi@VVV				\pi_0@VVV\\
U		@>\psi>>	\R^{n-k} \cross \R_+^k
\end{CD}
\end{equation}
where
$\pi_0 \colon (S^1\cross\R)^{n-k} \cross \C^k \to \R^{n-k} \cross \R_+^k$
is the projection to $\R$ on the first $n-k$ factors, and 
the projection $(x,y) \mapsto \frac{1}{2}(x^2 + y^2)$ 
in the last $k$ factors,
and where we take the standard symplectic structure on 
$\C$ and $S^1\cross\R \cong T^*S^1$.
\end{defnn}

\begin{rmk}
This implies that on an open subset of $B$, 
the preimages of points are Lagrangian manifolds.  
The level sets of $\pi$ form a singular Lagrangian fibration.
The singular fibres are those with $k>0$.
Note, however, that the fibres are only ``singular'' in terms of 
the fibration.  They are still smooth manifolds.
\end{rmk}

\begin{defn}
A \define{locally toric singular real polarization} 
on $M$ is the distribution $P \subset TM$
associated to a locally toric singular Lagrangian fibration 
(i.e., $P_x$ is the set of directions tangent to the leaf through~$x$).
\end{defn}

This definition includes toric manifolds, by the Local Normal Form.
In this case the map $\pi$ is the moment map, and $B$ can be 
taken to be the moment polytope.

It also includes more general 
integrable systems with certain kinds of singularities.
Eliasson in~\cite{eliassonthesis} and~\cite{eliassonelliptic}
and Miranda in~\cite{evathesis} established the local symplectic 
classification of non-degenerate singularities of integrable 
Hamiltonian systems: 
such singularities are isomorphic to the product of singularities of 
three basic types, called 
\emph{elliptic, hyperbolic,} and \emph{focus-focus.}
Definition~\ref{spacedef} includes integrable systems which have only elliptic
singularities, by the following theorem (due to Dufour and Molino and 
Eliasson, here taken from Zung~\cite{zung}):

\begin{thm}[3.9 in \cite{zung}]
Let $N$ be an elliptic singular leaf of codimension $k$
in an integrable system with moment map $F\colon M^{2n} \to \R^n$.
Then on a tubular \nbhd\ of $N$, there exist symplectic coordinates
$(x_1,\ldots,x_n,y_1,\ldots,y_n)$ so that:
\begin{itemize}
\item $y_1,\ldots,y_{n-k}$ are mod 1
\item $\w = \sum dx_j\wedge dy_j$
\item $N = \{ x_1 = \cdots = x_n = y_{n-k+1} = \cdots =y_n = 0 \} $,
i.e.\ $N$ is a $n-k$-torus, with coordinates $y_1,\ldots,y_{n-k}$
\item $F$ is a smooth function of $x_1,\ldots, x_{n-k}$, and 
$x_j^2 + y_j^2$ for $n-k < j \leq n$.
\end{itemize}
\end{thm}

(The relation of the notation in this theorem to that used in this paper 
is as follows: $k$ represents the same thing.
What we call $t_j, \theta_j$ in Section~\ref{s:cylinder} correspond to 
$x_j ,  y_j$ (except for factors of $2\pi$) for $j\leq n-k$.
What we call $s_j$ in Section~\ref{s:multdm}
is $\frac{1}{2}(x_i^2 + y_i^2)$ 
in the coordinates in this theorem, 
for $1\leq j \leq k$
, $i = (n-k)+j$.)

Kogan in \cite{Kogan} gives a description of the structure 
of these spaces.

\begin{defnn}\labell{Jdef}
Given a manifold $M$ with prequantization line bundle $\LL$ and 
(possibly singular) 
real polarization $P$,
a section $\sigma$ of $\LL$ over $U \subset M$
is \define{flat along the leaves,} or \define{leafwise flat}, if 
$ \nabla_X \sigma = 0 $
for all $X \in P$, at every point of $U$.

We denote the sheaf of leafwise flat sections on $M$ by $\J_M$ (or just $\J$).
\end{defnn}

\begin{defn}
The \define{quantization} of $M$ is the sum of the cohomology groups 
of $M$ with coefficients in the sheaf of leafwise flat sections:
\begin{equation}\labell{qn-def}
\Q(M) := \bigoplus_q H^q (M;\J) 
\end{equation}
\end{defn}

Thus, in this paper the central items of interest are 
the sheaf cohomology spaces $H^q(M;\J).$

\begin{rmk}
Various authors define the \qn\ in terms of sheaf cohomology, 
either as the direct sum as in~\eqref{qn-def} (e.g.\ in~\cite{JW}), 
or as the alternating sum of cohomology (e.g.\ in~\cite{ggk}).
In~\cite{GS} the authors call the groups 
$H^q(M;\J)$ ``the basic quantum data associated with $M$\ldots,'' 
without defining the \qn\ as either sum.
However, in all of these cases, as well as in this paper, all but one of 
the groups turn out to be zero, and so whether one takes the 
direct or the alternating sum doesn't matter in the end.
We use the convention of~\eqref{qn-def}, and call the resulting object 
``the quantization of $M$.''
\end{rmk}

\subsubsection{Bohr-Sommerfeld leaves and \'Sniatycki's theorem}\labell{ss:BS}

Let $M$ be a prequantizable compact symplectic manifold of dimension $2n$, 
as in the previous section.

\begin{defn}
A leaf $\lf$ of the polarization $P$ is a \define{\BS\ leaf} 
if there exists a globally defined nonzero section of $\LL$ along $\lf$, 
whose covariant derivative
(with respect to $\nabla$) is zero in directions tangent to $P$.
(Here ``globally defined'' means defined on all of $\lf$, 
not all of $M$.)

The \define{\BS\ set} is the set of points in $B$ whose preimages
are \BS\ leaves.
\end{defn}

In \cite{Sn}, \'Sniatycki proves that, in the case where the 
projection map $\pi \colon M \to B$ is a fibration,
the cohomology groups $H^q(M;\J)$
appearing in \eqref{qn-def} are all zero except 
in dimension $n$.  
Furthermore, $H^n$ can be computed by counting \BS\ leaves.  
More precisely, we have the following result:

\begin{thmn}[\'Sniatycki, 1975 \cite{Sn}]
Let $M$ be a $2n$-dimensional symplectic manifold, with a 
prequantization line bundle $\LL$ as above.  Let $P$ be a real polarization 
such that the projection map $\pi \colon M \to B$ is a fibration 
with compact fibres.  Then $H^q(M;\J) = 0$ for all $q\neq n$.

Furthermore, let 
$\Gamma_{BS}(\LL)$ be the space 
of smooth sections of $\LL$ along the union of \BS\ leaves, 
and $C_P^\infty(M)$ be the ring of 
functions on $M$ constant on leaves of $\pi$.
Then provided $P$ satisfies an orientability condition,
$H^n(M;\J)$ is isomorphic to $\Gamma_{BS}(\LL)$, 
as modules over $C_P^\infty(M)$.

More generally, if the leaves are not compact, then similar results
are true with $n$ replaced by the rank of the fundamental group 
of a typical integral manifold of $P$.
\end{thmn}

As a vector space, $H^n(M;\J)$ is isomorphic to the direct sum of 
copies of $\C$, with one copy for each \BS\ leaf.

\subsection{Examples}

\newcommand{\xfour}{(x_1,x_2,y_1,y_2)}

\begin{eg} A simple example of a compact space to which \'Sniatycki's 
results apply is $T^4$ fibred over $T^2$.  
If $T^4$ has coordinates $\xfour$, all mod $2\pi$, 
and standard symplectic form, then projecting to $T^2$ via 
\[ \xfour \mapsto (x_1,x_2) \]
is a Lagrangian fibration.  
\end{eg}

\vspace{1ex}

\begin{eg}
A less trivial example is Thurston's example \cite{Th}
of a symplectic manifold $M_\Theta$
which is not K\"ahler (and therefore not toric).  
Cannas da Silva in \cite{ana} gives the following description of $M_\Theta$:

Let $\Gamma$ be the discrete group generated by the following
symplectomorphisms of $\R^4$:
\begin{equation*}
\begin{split}
\gamma_1 &= \xfour \mapsto (x_1, x_2 + 1, y_1, y_2)\\
\gamma_2 &= \xfour \mapsto (x_1, x_2, y_1, y_2 + 1)\\
\gamma_3 &= \xfour \mapsto (x_1 + 1, x_2, y_1, y_2)\\
\gamma_4 &= \xfour \mapsto (x_1, x_2+ y_2, y_1 + 1, y_2)\\
\end{split}
\end{equation*}
Then $M_\Theta = \R^4 / \Gamma$, with
symplectic form  $\omega = dx_1\wedge dy_1 + dx_2 \wedge dy_2$.

If we map $M_\Theta$ to the 2-torus $T^2$ by 
\[ \xfour \mapsto (y_1,y_2) \]
(all coordinates taken mod 1) then the fibres are tori
in the $(x_1,x_2)$ coordinates, which are thus in fact Lagrangian 
submanifolds, so $M_\Theta$  is a 2-torus bundle over the 2-torus.
This fibration is a non-singular Lagrangian foliation, and so
$M_\Theta$ is a compact manifold that satisfies the hypotheses of 
\'Sniatycki's theorem.
\end{eg}

\vspace{1ex}

\begin{eg}
We can use Thurston's manifold to construct further 
(admittedly somewhat artificial) examples of manifolds which satisfy 
Definition~\ref{spacedef} but are not toric:
simply take the product  $M_T \cross M_\Theta$ 
of Thurston's manifold with any compact toric manifold $M_T$.
These do not have a global torus action of maximal dimension,
because of the $M_\Theta$ factor, but they are still locally toric, and so 
our results apply.
 \'Sniatycki's 
theorem does not apply, however, since the foliation is singular.
\end{eg}

\subsection{Aside: Rigidity of \BS\ leaves}

\begin{defn}
We say that $M$ satisfies \define{Bohr-Sommerfeld rigidity} 
if the \BS\ leaves of $M$ are independent of the  
choice of prequantum connection on $\LL$.
\end{defn}

\begin{prop}\label{rigid}
Let $i\colon \lf \into M$ be the inclusion of a leaf of the polarization.
If the induced map $i_* \colon H_1(\lf,\Z) \to H_1(M,\Z)$ is zero
for all leaves $\lf$, 
then $M$ satisfies \BS\ rigidity.
\end{prop}

\begin{proof}
If $i_*$ is zero, 
this means that any loop $\gamma$ on $\lf$ is homotopic to a point
in $M$.  

Let $\Sigma$ be a surface spanning $\gamma$.  
For definiteness, assume that the prequantization connection
has a potential one-form $\Theta$ defined on all of $\Sigma$.
Then from \eqref{holonomy}, 
the holonomy around $\gamma$ is given by 
\[ \text{hol}_\gamma = \exp \left( i \int_\gamma \Theta \right) \]
which equals 
\begin{equation}\labell{monodromy}
 \exp \left( i \int_{\Sigma} \omega \right) 
\end{equation}
by Stokes' theorem.
(If there is no $\Theta$ defined on all of $\Sigma$, break $\Sigma$
up into little surfaces over which $\Theta$ exists, apply Stokes' on 
each one, and piece back together.) 

Now $\lf$ is a Bohr-Sommerfeld leaf if and only if
$\text{hol}_\gamma=1$ for all loops $\gamma$ on $\lf$.  
By \eqref{monodromy}, this will be true if and only if
\[\exp\left( i \int_\Sigma \omega\right)\] 
is 1 for all $\Sigma$ which span a loop on $\lf$.
This depends only on $\omega$ and $\lf$, not on 
the connection form.  
\end{proof}

\begin{cor}
A compact symplectic toric manifold has rigid Bohr-Sommerfeld leaves.
\end{cor}

\begin{proof}
All odd-degree homology groups 
of a compact symplectic toric manifold are zero.
(See, for example, Theorem I.3.6 in \cite{ACL}.)
Thus the 
image of $H_1(\lf,\Z)$ in $H_1(M,\Z)$ is certainly zero, and so
$M$ satisfies Bohr-Sommerfeld rigidity.
\end{proof}

\begin{rmk}
Note that $\R\cross S^1$, considered in the next section,
 does not satisfy rigidity: if we 
change the connection by adding a (non-integer) 
constant multiple of $d\theta$ 
to it, this changes the Bohr-Sommerfeld leaves.
\end{rmk}

\section{The cylinder}\labell{s:cylinder}
The first model space we consider is the cylinder $\R\cross S^1$.  
In this section we compute its sheaf cohomology by hand.
For this section, let $M$ denote  
$\R\cross S^1$, with coordinates $(t,\theta)$, 
where $\theta$ is taken mod $2\pi$,
and symplectic form $\w = dt \wedge d\theta$.

In \ref{ss:BohrS}, we give the basic definition and set-up of the manifold,
describe the sheaf of sections flat along the leaves,
and calculate the \BS\ leaves.
In the next several sections, we find the sheaf cohomology of a 
simple type of subset: a ``band'' around the cylinder.  
Our strategy will be to find the cohomology of simple ``bands'' 
directly, and then piece them together using a Mayer-Vietoris sequence.

We start by computing the cohomology with respect to simple kinds
of covers (in~\ref{ss:sheaf-cylinder} and~\ref{s:advbrickwall}).
Next, in~\ref{ss:cyl-MV}, we derive a Mayer-Vietoris sequence for our
particular sheaf.  Finally, in~\ref{ss:cyl-refcov}, we take the limit over 
finer and finer covers to find the actual sheaf cohomology. 
(Cf.\ section~\ref{ss:sheaf-cohom}.)

\subsection{Flat sections and Bohr-Sommerfeld leaves}\labell{ss:BohrS}

$M$ is equipped with a  real \pol, given by vectors tangent to the $S^1$
directions. 
The leaves of the polarization are 
the fibres of the projection $M \to \R $.

Let $\LL = M\cross\C$ be the trivial bundle over $M$.
Let $\sigma$ be a section of $\LL$, which we can view as a 
complex-valued function.
Since $\w = dt \wedge d\theta = d(t\,d\theta)$,
a connection with potential one-form $t\,d\theta$
has curvature  $\w$.
Thus the connection given by
\begin{equation}\labell{cyl-conn}
\nabla_X \sigma =X(\sigma) - \sigma i t\, d\theta(X)
\end{equation}
makes $\LL$ into a prequantization line bundle over $M$.
(Note that, in this case, the potential one-form is defined on the 
entire manifold.)

In order to calculate the sheaf cohomology, we need to know 
which sections of $\LL$ are flat along the leaves 
(see Definition~\ref{Jdef}; as there, we denote the sheaf of 
such sections by $\J$).

\begin{prop}\labell{cyl-flat}
Let $U\subset M$, 
and let a section of $\LL$ over $U$ be given by 
a map $\sigma \colon U \to \C$.  
Then $\sigma$ is flat along the leaves if and only if 
it locally has the form
\begin{equation}\labell{flat}
 \sigma = a(t) e^{i t \theta}
\end{equation}
for some smooth function $a(t)$.
\end{prop}

\begin{proof}
This result follows directly from the 
 description of the connection given above in \eqref{cyl-conn}.
The section $\sigma$ is flat along the leaves if 
\[ 0 = \nabla_X \sigma = X(\sigma) - \sigma \, i t \,d\theta(X)
	\qquad \forall X \in P.  \]
The polarization $P$ is the span of $\frac{\del}{\del\theta}$,
so this condition is equivalent to 
\[
0= \nabla_{\!\ddtheta} \sigma = \frac{\del\sigma}{\del\theta}
	-\sigma i t,
\]
i.e.
\[
\frac{\del\sigma}{\del\theta} = i t \sigma.
\]
This is a differential equation for $\sigma(t,\theta)$, which 
is easily solved, giving \eqref{flat}.
\end{proof}

Recall from Section~\ref{ss:BS} that 
a leaf $\lf$ satisfies the \emph{Bohr-Sommerfeld condition}
if it possesses a global covariant constant section (other than zero),
i.e. a section defined on the whole leaf which is flat along the leaf.

\begin{lemma}\labell{l:cyl-BS}
The Bohr-Sommerfeld set is $\Z \subset \R$,
and the Bohr-Sommerfeld leaves are $\{m\}\cross S^1$, $\; m \in \Z$.
\end{lemma}

\begin{proof}
Fix a leaf $\lf_0 = \{t_0\}\cross S^1$.  
By \eqref{flat}, a flat section $\sigma$ over $\lf_0$
is of the form $c e^{i t_0 \theta}$ for some constant $c$.  
It will always exist on a small \nbhd\ in the leaf, 
but will only be defined on the whole leaf if its values are 
compatible as it wraps around the leaf.  That is, values of $\sigma$ must 
agree for values of $\theta$ that differ by $2\pi$.  
This requires $e^{it_0\theta} = e^{it_0(\theta+2\pi)}$, i.e.\
$e^{2\pi i t_0} = 1$, i.e.\ $t_0 \in \Z$.
\end{proof}

Note that the space of global covariant constant sections over one leaf
is one-dimensional: $\{ \sigma = c e^{it_0\theta} \mid c\in \C \}$.

Let $I\subset\R$ be an open interval, and let $U = I\cross S^1 \subset M$.
By \'Sniatycki's theorem, 
\[ 
H^1(U,\J) \cong \bigoplus_{m \in \Z\cap I} \C; 
\qquad H^k(U,\J) = 0,\; k\neq 1.
\]

In the following sections we compute the sheaf cohomology of $U$ 
directly, and show that it agrees with \'Sniatycki's theorem.

\subsection{Sheaf cohomology}\labell{ss:sheaf-cylinder}
In this section we compute the \v Cech cohomology, 
with coefficients in $\J$, 
of a band in $M$, 
with respect to a particular cover.  
As a warm-up, in this section we use the simplest
possible cover; however, as we will see laer,
this case already shows all the important 
features of the calculation.

\begin{defn}
A \define{band} is a subset of $M$ of the form $I \cross S^1$, 
with $I\subset \R$ a bounded open interval.
\end{defn}

\begin{figure}[htbp]
\centerline{\includegraphics[width=2in,height=1in]{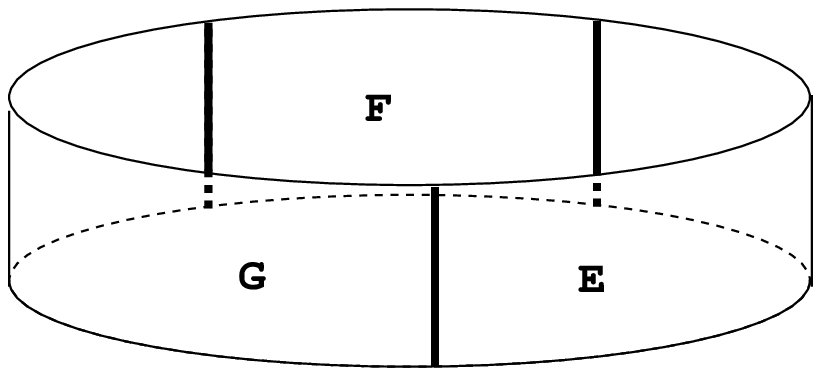}}
\caption{The cover $\E_3$}\label{E3}
\end{figure}

Let $U$ be a band around the cylinder that 
contains at most one Bohr-Sommerfeld leaf.
Partition $U$ into three
rectangles  $E$, $F$, and $G$ by partitioning $S^1$ 
into three intervals.   (See Figure~\ref{E3}, where the heavy line
indicates an overlap.)
We will calculate the 
cohomology of $U$ with respect to the cover $\E_3 = \{ E,F,G \}$.

\subsubsection{$H^0$}

We can see directly that $H^0(U;\J) = 0$, since $H^0(U;\J)$ is just 
the space of global sections of $\J$ over $U$, which we know 
from the argument in the proof of Lemma~\ref{l:cyl-BS} is $\{ 0\}$.  
However, we will calculate it directly to begin seeing how the 
\v Cech approach works in this situation.

A \v{C}ech 0-cochain $\eta$ is an assignment, to each of the sets $E$, $F$,
and $G$, of a flat section over that set.  Such a section will have the 
form $a(t) e^{it\theta}$.  Index the section by the set, so 
the piece of $\eta$ on $E$ is $a_E(t) e^{it\theta}$, etc.

On each set, the coordinate $\theta$ can be defined, 
even though can not be defined on all of $S^1$.  Fix a branch of 
$\theta$ on each open set $W$, and denote it by $\theta_W$.
Choose these branches so that $\theta_F = \theta_E$ on $E\cap F$, 
$\theta_G = \theta_F$ on $F\cap G$, and 
$\theta_G = \theta_E + 2\pi$ on $G\cap E$.

The coboundary of $\eta$ is the collection 
\[ (\delta \eta)_{VW} = \eta_W - \eta_V 
= a_W(t)\, e^{it\theta_W} - a_V(t)\, e^{it\theta_V}, \]
and so $\eta$ will be a cocycle if each of these are zero.  
Applying this to each of the three sets, we have
that $\eta$ is a cocycle if{f}
\begin{equation}\labell{0-coc}
\begin{split}
0 &= a_F(t)\, e^{it\theta_F} - a_E(t)\, e^{it\theta_E}
	\qquad\text{on } E\cap F\\
0 &= a_G(t)\, e^{it\theta_G} - a_F(t)\, e^{it\theta_F}
	\qquad\text{on } F\cap G\\
0 &= a_E(t)\, e^{it\theta_E} - a_G(t)\, e^{it\theta_G}
	\qquad\text{on } G\cap E
\end{split}
\end{equation}
In the first two of these equations, the $\theta$ coordinates are
equal on the intersections, and so we can cancel the exponential 
factors; this implies 
\begin{align}\labell{0coc1}
a_E(t) = a_F(t) \qquad\text{and}\qquad a_F(t) = a_G(t).
\end{align}
However, $\theta_E$ and $\theta_G$ differ by $2\pi$ on the intersection,
and so the third equation becomes 
\[ 0 = a_E(t)\, e^{it\theta_E} - a_G(t)\, e^{it\theta_E+2\pi it}, \]
which implies 
\begin{equation}\labell{0coc2}
a_E(t) = e^{2\pi it} a_G(t).
\end{equation}
Conditions \eqref{0coc1} and \eqref{0coc2} together 
require that $e^{2\pi it}=1$,
which cannot happen on an interval of $t$-values.  
Thus there are no 0-cocycles, and  $H^0 = 0$.

\subsubsection{$H^1$}

The one-dimensional case is more interesting.  A 1-cochain $\eta$
is an assignment of a flat section $\eta_{VW}$ to each 
intersection $V\cap W$; thus a 1-cochain is given by a triple
of functions
\[
\{ a_{EF}(t)e^{it\theta}, a_{FG}(t)e^{it\theta}, a_{GE}(t)e^{it\theta} \}. 
\]
Since there are no triple intersections in this cover, there
are no 2-cochains, and thus every 1-cochain is a cocycle.  

Note that, since $\eta$ is determined by the $a$'s, and each $a$ 
is a smooth function of $t$ on $I$, the space of cocycles is isomorphic
to $C^\infty(I)^3$.   

We now consider when a 1-cochain is a coboundary, 
namely, when there exists a 0-cochain 
$\beta = \{ b_E e^{it\theta_E},  b_F e^{it\theta_F},  
b_G e^{it\theta_G} \}$ with 
$\delta \beta = \eta$.  
This requires that 
\begin{equation}\labell{eq:beta}
\eta_{VW} = \beta_W - \beta_V \quad \text{on } V\cap W,
\end{equation}
for each pair of $V$ and $W$.

\noindent {\bf Notation:}
We will write $EF$ for $E\cap F$, and so on.  
The order in which intersections are written matters in
\v Cech cohomology, for bookkeeping, and so we need to be consistent.  
We will write intersections in the order $EF$, $FG$, $GE$.
Furthermore, we will 
use the convention that, on any intersection $VW$, we use the $\theta$ 
coordinate from $V$.

Applying these conventions 
to the possible equations~\eqref{eq:beta}, and using the formulas for 
$\eta$ and $\beta$, we obtain that $\delta \beta = \eta$ if{f} 
\begin{subequations}\labell{1-coc}
\begin{align}
a_{EF}(t) e^{it\theta_E} &= b_F(t) e^{it\theta_F} - b_E(t) e^{it\theta_E}
\qquad \text{on }EF\\
a_{FG}(t) e^{it\theta_F} &= b_G(t) e^{it\theta_G} - b_F(t) e^{it\theta_F}
\qquad \text{on }FG\\
a_{GE}(t) e^{it\theta_G} &= b_E(t) e^{it\theta_E} - 
	b_G(t) e^{it\theta_G}\labell{c3}
\qquad \text{on }GE
\end{align}
\end{subequations}
In each of these equations, all the $\theta$ coordinates are equal
on the relevant set, except in \eqref{c3}, where they 
differ by a factor of $2\pi$.  
Applying this fact, and cancelling common factors of $e^{it\theta}$, we obtain
the following system of three equations
\begin{equation}\labell{coc-system}
\begin{split}
a_{EF} &= b_F - b_E\\
a_{FG} &= b_G - b_F\\
a_{GE} &= e^{-2\pi i t} b_E - b_G
\end{split}
\end{equation}
in the three unknown functions $b_E$, $b_F$, and $b_G$ on $I$.
They must be true for each value of $t$ in $I$.
The matrix of this system is 
\begin{equation}\labell{matrix}
 \begin{bmatrix}
-1&  1&  0\\
0&  -1&  1\\
e^{-2\pi it}&  0& -1
\end{bmatrix} 
\end{equation}
which is invertible precisely when $e^{-2\pi it} \neq 1$.
Thus, by elementary linear algebra, 
the system \eqref{coc-system} has a solution, for any collection
of the $a_{VW}$, if $e^{-2\pi it}$ is never 1 on $U$.
In this case, every cocycle is a coboundary, and $U$ has trivial cohomology.

The other possibility is if $e^{2\pi it}=1$ somewhere in $I$,
which only occurs if $I$ contains an integer $m$.
In this case, by further linear algebra, the system \eqref{coc-system} 
only has a solution if $\eta$ satisfies the condition
\begin{equation}\labell{cob-cond}
a_{EF}(m) + a_{FG}(m) + a_{GE}(m) = 0.  
\end{equation}
Thus $\eta$ is a coboundary precisely when \eqref{cob-cond} holds, and 
so the cohomology of $U$ is
\begin{equation}\labell{cyl-H1}
H^1 = C^\infty (I)^3 /
   \{a_{EF}(m) + a_{FG}(m) + a_{GE}(m) = 0 \}. 
\end{equation}

\begin{lemma}\labell{coh-quot}
The quotient \eqref{cyl-H1} is isomorphic to $\C$.
\end{lemma}
\begin{proof}
Map $C^\infty(I)^3 \to \C$ by $\phi(f,g,h) = f(m) + g(m) + h(m)$,
where $m$ is the (unique) integer in $I$.  
This surjective 
homomorphism gives the desired isomorphism.
\end{proof}


Finally, note that the condition $e^{2\pi it}=1$ occurs precisely 
at the Bohr-Sommerfeld leaves.  Therefore, we have proved the 
first part of the following result.
The second statement follows since 
for  $k > 1$, there are no $(k+1)$-fold intersections in 
this cover.

\begin{prop}\labell{e3-band-cohom}
Let $U$ be a band around the cylinder.  Then the sheaf cohomology 
 of $U$ with respect to the cover $\E_3 = \{ E,F,G \}$ is
trivial if $U$ does not contain a \BS\ leaf.  
If $U$ contains one \BS\ leaf, its cohomology is 
\begin{equation*}
H^k_{\E_3} (U;\J) \cong 
\begin{cases}
  \C &k = 1\\
  0   &k\neq 1
\end{cases}
\end{equation*}
\end{prop}

This is precisely what we expect from \'Sniatycki's results.

\begin{lemma}\labell{Ek}
Let $\E_k$ be a cover of the band similar to $\E_3$, but with $k$ 
sets instead of 3.  Then the cohomology calculated with respect to $\E_k$ 
is the same as that calculated with respect to $\E_3$.
\end{lemma}

\begin{proof}
In this case, the same argument as for the cover $\E_3$ applies;
the only difference is that we have $k$ equations instead of 3 
in \eqref{0-coc} and \eqref{1-coc}, and the matrix \eqref{matrix}
is a $k\cross k$ matrix.  Its determinant is the same, however, (up to sign) 
and we obtain the same conclusion as in 
Proposition~\ref{e3-band-cohom} for cohomology 
with respect to the cover $\E_k$.
\end{proof}

\subsection{Brick wall covers}\labell{s:advbrickwall}
Eventually (see \ref{ss:cyl-refcov} below), 
we will find the cohomology of a band $U$ by breaking 
it up into ``sub-bands,'' finding the cohomology of each one by hand, 
and then piecing them together using Mayer-Vietoris.  
If each band has a cover of the form $\E_k$, the cover induced on their
intersection will be similar to $\E_k$ but with two ``layers,'' 
and thus the calculation of section~\ref{ss:sheaf-cylinder} 
is not sufficient.  
In this section, we define a type of covering we call a ``brick wall,'' 
which looks like $\E_k$ but with more layers, 
and compute 
the cohomology of a band with this type of cover.  

Later (in~\ref{ss:cyl-refcov}), 
we will take the direct limit over all covers of $U$ 
(see~\ref{ss:sheaf-cohom}) to find the actual sheaf cohomology.
We will use Lemma~\ref{cofinal} and show that the set of 
brick wall covers is cofinal in the set of all covers of $U$ 
(see Lemma~\ref{brickwall} below),
and so it will be sufficient for all our purposes to consider only 
brick wall covers.

It will be enough to consider covers with only two layers of bricks,
since the cohomology of covers with more layers of bricks can be 
found by piecing together two layers at a time, using Mayer-Vietoris.

\begin{defn}
A \define{brick wall} cover of a band in the cylinder 
(or, more generally, of any rectangle)
is a finite covering by open $t$-$\theta$ rectangles (``bricks''),
satisfying the following properties: 
\begin{itemize}
\item The rectangles can be partitioned into sets (``layers'') so that 
all rectangles in one set cover the same interval of $t$ values
(``All bricks in the same layer have the same height'');
\item Each brick contains points that are not in any other brick; and
\item There are no worse than triple intersections, i.e., the 
intersection of two bricks in one layer does not meet the intersection
of two bricks in either of the two adjoining layers.
\end{itemize}
\begin{figure}[htbp]
\centerline{\includegraphics[width=1.5in,height=1in]{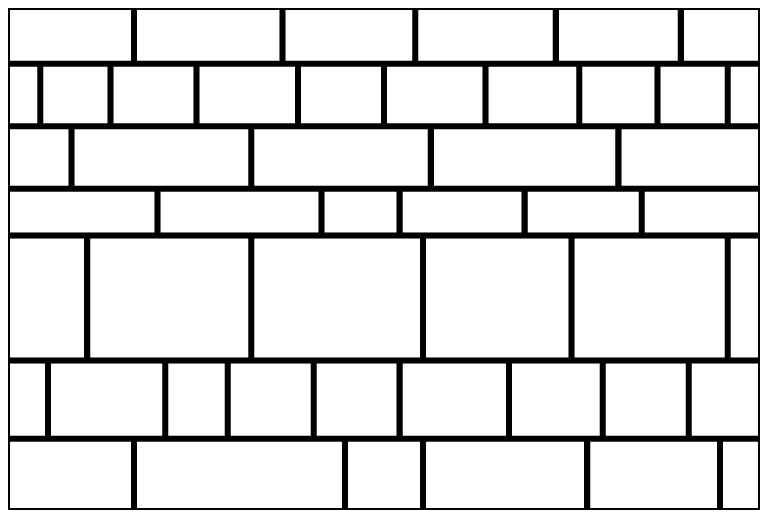}
\qquad \qquad \includegraphics[width=1.5in,height=1in]{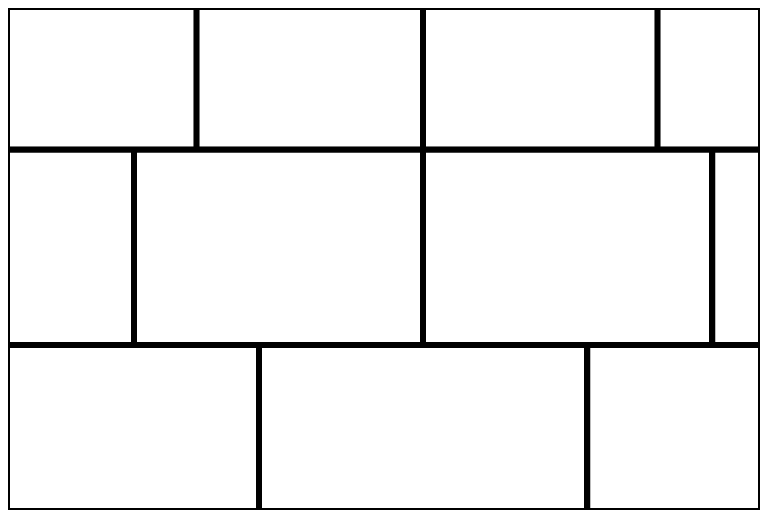}}
\caption{A brick wall cover, and one which is not}\labell{brickfig}
\end{figure}
Note that we do not require that the number of bricks be the same in each
layer, nor that the layers have the same height, nor that the bricks within
one layer have the same width.
See Figure~\ref{brickfig}, where we have ``unrolled'' the band, 
and where thick lines indicate intersections.
Later we will allow brick walls with countably many layers, although 
the number of bricks in each layer will still be finite.
\end{defn}

Let $U=I\cross S^1$ be a band around the cylinder, which we cover
by a brick wall of two layers.  
Let the top layer have $m$ bricks $A_1$ through $A_m$, 
and the bottom layer have $n$ bricks $B_1$ through $B_n$.
(We choose our numbering of the $A$'s and $B$'s so that $B_1$ meets both
$A_1$ and $A_m$.)
Denote this covering by $\B^m_n$.
Let $I_A$ and $I_B$ denote the intervals of $t$ values which are covered 
by the $A$ and $B$ layers respectively.  
We also assume that $U$ contains at most one \BS\ leaf, and that this 
leaf is not contained in the intersection between the two layers.
(It is possible to perform the calculation without this assumption. 
However, this case is sufficient for our purposes---since eventually we'll
be taking finer and finer covers, we can always arrange that the \BS\ 
leaves avoid the intersections between brick layers---and avoids 
some complications in the argument.)  
The most challenging part of this calculation is the bookkeeping, 
so we will set out our notational conventions at the outset.  

A 0-cochain is given by a leafwise flat section on each brick in the
cover, which as we showed in \ref{ss:BohrS} is of the form $f(t) e^{it\theta}$
for some smooth function $f$.  We will denote the functions
corresponding to the set $A_j$ and $B_k$ by $a_j$ and $b_k$,
respectively, so that for example $a_k(t) e^{it\theta}$ is the element
of the cochain defined on $A_k$.

A 1-cochain is given by a section on each intersection of two bricks.  
The functions corresponding to the intersections $A_i \cap A_j$ 
and $B_k \cap B_l$ 
will be denoted by $a_{ij}$ and $b_{kl}$ respectively.
The function corresponding to the set $A_j \cap B_k$ will be denoted 
by $c_{jk}$.  
Thus, the $a_{ij}$'s and $b_{kl}$'s give the parts of the 1-cochain
defined on the intersections \emph{within} one layer of the brick wall, 
while the $c_{jk}$'s are on the intersections \emph{between} the layers.
In particular, the $a$'s are functions of $t$ defined on $I_A$, the $b$'s are 
defined on $I_B$, and the $c$'s are defined on $I_{AB} = I_A \cap I_B$.

As we noted before, 
the order in which sets are written in intersections is important
in \v Cech cohomology---it doesn't matter 
how we do it, as long as we're consistent---and 
so we set out our conventions here.  
We will write sets with smaller index before those with larger index
(wrapping around, so that $n$ is considered ``less'' than 1), 
and sets on the $A$ layers will be written before sets in the $B$ layer.
(Thus, we would write $A_2A_3$, $A_mA_1$, and $A_mA_1B_1$.)
Also, on an intersection, we will use the $\theta$ coordinate from the 
set written \emph{first} in the intersection by these conventions.
Finally, the ``branches'' of $\theta$ will be chosen so that
$\theta_{A_m} = \theta_{A_1} + 2\pi$ on $A_mA_1$, 
$\theta_{B_n} = \theta_{B_1} + 2\pi$ on $B_n B_1$, 
$\theta_{B_1} = \theta_{A_1}$ on $A_1 B_1$, and the $\theta$ coordinates 
on all other double intersections agree.

\gap

As in section \ref{ss:sheaf-cylinder}, $H^0 = 0$ because there are no 
global sections, and $H^j = 0$ for $j \geq 3$ because this cover has 
at most triple intersections.  
We start with $H^1$---as 
with the simpler cover, this is where all the interesting things 
happen---and deal with $H^2$ later.

Suppose we're given a 1-cocycle.
This is a collection of $a_{ij}(t)$, $b_{kl}(t)$, and $c_{jk}(t)$, 
for all possible intersections of the $A$'s and $B$'s,\footnote{
Note that the $a_{ij}$ and $b_{kl}$ will actually be 
\[ a_{i(i+1)} \quad\text{and}\quad b_{k(k+1)}, \qquad 1\leq i \leq m-1, 
\quad 1 \leq k \leq n-1, \]
plus $a_{m1}$ and $b_{n1}$.
The $c_{jk}$, on the other hand, will be defined for all pairs 
$(j,k)$ for which $A_j$ and $B_k$ intersect; it is not possible 
a priori to say which pairs exist, but as discussed in a moment,
there will be $m+n$ of the $c_{jk}$.} 
satisfying certain conditions, which we'll deal with in a moment.
We seek $a_j$, $b_k$ defined on each $A_j$ and $B_k$ which make up 
a 0-cochain whose coboundary is our $\{ a_{ij}, b_{kl}, c_{jk} \}$.

Just looking at the bricks within one layer, the situation is identical
to the cover $\E_k$ considered in the previous section, 
and we get a system of equations like \eqref{1-coc}.
For the $A$ layer, we get $m$ equations in the $m$ unknown functions  
$a_1, \ldots, a_m$ defined on $I_A$, which, after
applying the relationships between the various $\theta$ coordinates and 
cancelling common factors of $e^{it\theta}$, 
give the $m$ equations
\begin{equation}\labell{aij}
\begin{split}
a_{12} &= a_2 - a_1 \\
a_{23} &= a_3 - a_2 \\
\vdots \\
a_{m1} &= a_m - e^{2\pi it} a_1 
\end{split}
\end{equation}
as functions of $t$ defined on $I_A$.
Provided these equations are consistent, 
they uniquely determine $a_1, \ldots, a_m$ from the $a_{ij}$
(by the same linear algebra argument as for the cover $\E_k$ 
in Lemma~\ref{Ek}).
As in in that case, these equations will be consistent
provided $e^{2\pi it} \neq 1$ on $I_A$.

A similar set of $n$ equations:
\begin{equation}\labell{bij}
\begin{split}
b_{12} &= b_2 - b_1 \\
\vdots \\
b_{n1} &= b_n - e^{2\pi it} b_1 
\end{split}
\end{equation}
holds for the $n$ functions $b_k$ on $I_B$,
which, again provided $e^{2\pi it} \neq 1$ (on $I_B$), 
uniquely determine the $b_k$ from the $b_{kl}$.

Thus, all of the functions making up our 0-cochain
are already determined just from the elements of the 1-cocycle 
which only exist within one layer.
However, we also have a number of equations coming from the intersections
\emph{between} the layers, 
which need to be taken into account. 

\begin{figure}[htbp]
\centerline{\includegraphics[width=3in,height=0.5in]{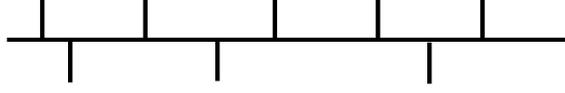}}
\caption{The intersection of two brick layers}\label{brickbranch}
\end{figure}

First of all, note that there are $n+m$ double intersections between
the layers.  This can be seen easily from Figure~\ref{brickbranch}, 
which shows the view of the middle of a typical two-layer 
brick wall cover.  Start at one point and go around the cylinder,
counting double intersections.  A new one will be added to the count
every time we pass one of the vertical lines, i.e.\ the intersection
between two bricks in the same layer.  Since there are $m$ bricks in 
one layer and $n$ in the other, there are a total of $m+n$ vertical
lines, and thus $m+n$ double intersections between bricks in 
different layers.  
(A similar argument shows that there are $m+n$ triple intersections
in this cover.)

On each such double intersection, there 
is a $c_{jk}$ which must satisfy
\begin{equation}\labell{cab} 
c_{jk} = b_k -a_j 
\end{equation}
as functions of $t$ on $I_{AB}$,
and so we have $m+n$ $c_{jk}$'s and $m+n$ such equations.

We also have relations among the $a_{ij}$, $b_{ij}$, and $c_{jk}$, 
coming from the fact that they make up a cocycle.  
These equations come from the triple intersections of sets in the cover.
Since there are $m+n$ triple intersections,
we have $m+n$ equations, which are of the form 
\begin{equation}\labell{cij}
\begin{split}
c_{(i+1)\, k} - c_{ik} + a_{i\, (i+1)} &= 0 
\qquad \text{from } A_i A_{i+1} B_{k}\\
\text{or}\\
b_{k\, (k+1)} - c_{ik} + c_{i\, (k+1)} &= 0 
\qquad \text{from } A_i B_{k} B_{k+1}
\end{split}
\end{equation}
depending on whether the intersection has two $A$'s or two $B$'s.  
(This is ignoring, for the moment, when there are factors of $e^{2\pi it}$ 
to worry about, which only happens near $A_1$ and $B_1$.)
These are the cocycle conditions, which are again equalities of functions 
of $t$ defined on $I_{AB}$.

Essentially, we have $m+n$ extra equations that the $a_i$ and $b_j$ 
must satisfy, but we have $m+n$ conditions that the extra equations 
satisfy, which is enough to cancel each other out.
In what follows we show more detail, but this is the essential idea.

Re-write equations \eqref{cij} as
\begin{subequations}\labell{cij2}
\begin{align}
c_{(i+1)\, k} - c_{ik} &= - a_{i\, (i+1)}\labell{cij-1} \\
c_{ik} - c_{i\, (k+1)} &= b_{k\, (k+1)} \labell{cij-2}
\end{align}
\end{subequations}
Think of the $a_{ij}$ and $b_{kl}$ as being given, and these $m+n$ 
equations as defining the $c_{jk}$ in terms of them.
Since there are $m+n$ unknowns (the $c_{jk}$), provided they are consistent 
and not underdetermined, they define the $c_{jk}$ uniquely.  
This will be the case provided the determinant of the matrix of coefficients 
of the system is not zero.  

Order the $c_{jk}$ by just going around the central part of the band.
From Figure~\ref{brickbranch}, we can see that each double intersection
meets two triple intersections; thus, each $c_{jk}$ appears in two 
(successive) 
equations in \eqref{cij2}: once with its predecessor, and once with 
its successor.  
If we write all of the equations in the form given in \eqref{cij2}, 
then a given $c_{jk}$ has a positive sign when it appears 
with its predecessor, and a negative sign when it appears with its 
successor, and these are its only appearances.
Thus the coefficient matrix for the system \eqref{cij2} will have 
mostly zeros, except each row will have one $-1$ and one $+1$ in 
adjacent entries.

Finally, there are two places  where we need to take into account the 
change in $\theta$ coordinates, which happens between $A_m$ and $A_1$ 
and between $B_n$ and $B_1$.  The two corresponding cocycle equations are
\begin{subequations}\labell{cocoord}
\begin{align}
c_{11} e^{it\theta_{A_1}} - c_{m1} e^{it\theta_{A_m}} 
	+ a_{m1}e^{it\theta_{A_m}} &=0 
	\qquad\text{on } A_m A_1 B_1 \labell{coco-1}\\
b_{n1} e^{it\theta_{B_n}} - c_{k1} e^{it\theta_{A_k}} 
	+ c_{kn}e^{it\theta_{A_k}} &=0 
	\qquad\text{on } A_k B_n B_1 \labell{coco-2}
\end{align}
\end{subequations}
(where $A_k$ is the brick in the $A$ layer that straddles the intersection
of $B_1$ and $B_n$).

Using the fact that $\theta_{A_1} = \theta_{A_m} - 2\pi$, we rewrite 
\eqref{coco-1} as 
\begin{equation*}
c_{11} e^{it\theta_{A_m}-2\pi i t} - c_{m1} e^{it\theta_{A_m}} 
   + a_{m1}e^{it\theta_{A_m}}=0  \qquad\text{on } A_m A_1 B_1 \\
\end{equation*}
which yields
\begin{equation}
c_{11} e^{-2\pi it} - c_{m1} = -a_{m1}.
\end{equation}
As for \eqref{coco-2}, since $\theta_{B_n} = \theta_{A_k}$, we can cancel 
the $e^{it\theta}$ terms immediately to get 
\begin{equation}
c_{kn} - c_{k1} = b_{n1}.
\end{equation}

Therefore, the coefficient matrix for the system \eqref{cij2} can be 
put in the form:
\begin{equation}\labell{cijmatrix}
\begin{bmatrix}
-1&  1&  0&  \cdots& 0& 0\\
0&  -1&  1&  \cdots& 0& 0\\
\vdots& & & \ddots\\
0&   0&  0&  \cdots& -1& 1\\ 
e^{-2\pi it}&  0& 0& \cdots& 0& -1
\end{bmatrix}
\end{equation}
By expanding along the bottom row, we can see that this matrix has 
determinant  $e^{-2\pi it} -1$, and thus the equations \eqref{cij} 
have a unique solution for the $c_{jk}$ provided, as usual, 
that $e^{2\pi it} \neq 1$ on $I_{AB}$---which is precisely what we 
are assuming.

Finally, it is straightforward to check that $c_{jk} = b_k - a_j$,
where the $a_j$ and $b_k$ are the ones found already, 
gives a solution to \eqref{cij}.

The upshot of all of this is that the $a_i$ and $b_k$ are 
determined entirely by the parts of the cocycle defined 
on the intersections within one layer,
namely the $a_{ij}$ and the $b_{kl}$.
The parts of the cocycle defined on the intersections 
between the layers (the $c_{jk}$) don't have any effect 
on the $a_i$ and $b_k$, because of the cocycle conditions.
Thus, given a 1-cocycle, provided that $e^{2\pi it}$ is never 1 on 
the band, it is the coboundary of a 0-cochain,  
and so the first \v Cech cohomology is zero.

Now consider the case when  $e^{2\pi it} = 1$ somewhere on the band.
If $e^{2\pi it_0} = 1$ for some $t_0 \in I_A$, then 
by an elementary linear algebra argument, the equations
\eqref{aij} are only consistent if 
\[ a_{12}(t_0) + a_{23}(t_0) + \cdots + a_{m1}(t_0) = 0. \]
In this case, there is a unique solution for $\{ b_1, \ldots, b_n \}$ 
on $I_B$, but there is only a solution for $\{ a_1, \ldots, a_m \}$ 
if $a_1(t_0) + \cdots + a_m(t_0) =0$.
Since $m \notin I_{AB}$, the system \eqref{cij2} still has a unique 
solution, which is compatible with the solutions for the $a$'s and $b$'s
by the cocycle conditions.  Thus, in this case, $H^1$ will be given by 
\begin{equation}\labell{quot1}
H^1 = Z^1 / \{ a_1(t_0) + \cdots + a_m(t_0) =0 \}
\end{equation}
where $Z^1$ is the space of 1-cocycles.

\begin{lemma}\labell{Z1}
The set of 1-cocycles $Z^1$ is isomorphic to 
\[ C^\infty(I_A)^m \oplus C^\infty(I_B)^n, \] 
provided there is no integer in $I_{AB}$.
\end{lemma}

\begin{proof}
As discussed above, if we think of the $\{ a_{ij} \}$ and $\{ b_{kl} \}$
as being given, we can view the equations \eqref{cij} as defining the
functions $\{ c_{jk} \}$ in terms of the $a$'s and $b$'s; as noted 
above, this system will have a unique solution if $e^{2\pi it} \neq 1$ 
on $I_{AB}$.  Thus, specifying a 1-cocycle amounts to giving the 
$m$ functions $\{ a_{12}, a_{23}, \ldots, a_{m1} \}$ on $I_A$ 
and the $n$ functions $\{ b_{12}, b_{23}, \ldots, b_{n1} \}$ on $I_B$.
Thus  $Z^1 \cong C^\infty(I_A)^m \oplus C^\infty(I_B)^n$.
\end{proof}

\begin{lemma}
The quotient \eqref{quot1} is isomorphic to $\C$.
\end{lemma}

\begin{proof}
Map  $Z^1 \cong C^\infty(I_A)^m \oplus C^\infty(I_B)^n \to \C$ via 
\[ (f_1, \ldots f_{m+n}) \mapsto f_1(t_0) + f_2(t_0) + \cdots f_m(t_0). \]
As in the argument in Lemma~\ref{coh-quot} 
(Section~\ref{ss:sheaf-cylinder}),
this is a surjective homomorphism, and gives the desired quotient.
\end{proof}

On the other hand, if $e^{2\pi it} = 1$ on $I_B$, 
the same argument applies to the $b_{kl}$ and gives the same result for
the cohomology. 

\gap

To sum up, we have shown the $k=0$ and $k=1$ cases of the following:

\begin{lemma}
Let $U$ be a band around the cylinder containing at most one \BS\ leaf,
with a brick wall cover $\B^m_n$ of two layers.
Assume that the \BS\ leaf is contained in at most one layer of bricks.
The the sheaf cohomology of $U$ with respect to the cover $\B^m_n$ 
is 0 if $U$ contains no \BS\ leaf, and 
\begin{equation}
H^k_{\B^m_n} (U;\J) \cong 
\begin{cases}
  \C &k = 1\\
  0   &k\neq 1
\end{cases}
\end{equation}
if it contains one \BS\ leaf.
\end{lemma}

\begin{proof}[Proof of $H^{\geq 2}$]
Unlike the simpler cover of section~\ref{ss:sheaf-cylinder},
we do have triple intersections in this cover, and so it is not 
immediate that $H^2 = 0$.  However, it will not take long to dispose of 
this calculation.

A 2-cochain is one section for each triple intersection, which 
again can be represented by a smooth function of $t$ on the appropriate
interval.  Since all the triple intersections lie along the intersection
of the two layers, these functions will be defined on $I_{AB}$.
Denote these functions by $f$'s, and call them $f_{j\,(j+1)\,k}$ 
for the intersections $A_j A_{j+1} B_k$ and $f_{j k\, (k+1)}$ 
for $A_j B_k B_{k+1}$.

Since there are no quadruple intersections, there are no 3-cochains, 
and so every 2-cochain is a cocycle.  Thus, to show that $H^2$ is trivial,
we need to show that every 2-cochain is a coboundary.  
So suppose we are given the $f$'s making up a 2-cochain.

There are $m+n$ triple intersections, as noted above, and thus $m+n$ $f$'s.  
On each intersection we have an equation of the form 
\begin{equation}
c_{(j+1)\, k} - c_{jk} + a_{j\, (j+1)} = f_{j\,(j+1)\,k} 
\end{equation}
or 
\begin{equation}
b_{k\,(k+1)} - c_{jk} + c_{j\,(k+1)} = f_{j\,k\,(k+1)}
\end{equation}
This gives $m+n$ equations in the $2(m+n)$ unknowns 
(the $a$'s, $b$'s, and $c$'s)
and so has infinitely many solutions. 
(The equations are clearly linearly independent and so consistent.)
Thus every 2-cocycle is a coboundary, and $H^2= 0$.

Finally, $H^k = 0$ for $k\geq 3$ since there are no $k+1$-fold intersections.

\end{proof}

\subsection{Mayer-Vietoris}\labell{ss:cyl-MV}

\newcommand{\ba}{{\boldsymbol{\alpha}}}

In this section we apply a Mayer-Vietoris type argument to find 
the cohomology of the union of two open sets.  The argument follows 
that given in \cite{BT}, Section I.2.  

Let $U=S^1 \cross I_U$ and $V=S^1\cross I_V$ be bands around $M$, 
where $I_U$ and $I_V$ are intervals.  
Fix a covering $\mathcal{A} = \{A_\alpha\}$ of $M$.  This induces coverings 
on $U$, $V$, and $U\cap V$ 
(which we will still denote by $\A$).
In what follows, we compute sheaf
cohomology always with respect to these covers.

Consider the sequence of \v{C}ech cochain complexes 
\begin{equation}\labell{ses}
0 \to C^*_\A (U \cup V, \J) \stackrel{q}{\to}
C^*_\A (U,\J) \oplus C^*_\A (V,\J) \stackrel{r}{\to}
C^*_\A (U \cap V, \J) \to 0
\end{equation}
For ease of notation, use $\ba$ for $\alpha_0\cdots\alpha_k$, so 
that $A_\ba = A_{\alpha_0\cdots\alpha_k}$. 
The map $q$ restricts a component of a cochain
$f_\ba$ on $A_\ba$ to $U\cap A_\ba$ and $V\cap A_\ba$, 
respectively.  
The map $r$ is defined as follows: 
\[ \bigl( r(f,g)\bigr)_\ba = f_\ba - g_\ba
\quad \text{on } U\cap V\cap A_\ba. 
\]

\begin{prop}\labell{p:exact}
The sequence \eqref{ses} is exact.
\end{prop}

\begin{proof}
By construction $rq=0$, so we wish to show $\ker r \subset$ im$\,q$.
Suppose $(f,g) \stackrel{r}{\mapsto} 0$.
This means that for each $\ba$, 
 $f_\ba = g_\ba$ on $A_\ba \cap U \cap V$.
By the first axiom of sheaves (see~\ref{ss:sheaf-cohom}), 
there exists a function 
$h_\ba$ on $A_\ba \cap (U\cup V)$ which restricts to $f_\ba$ 
and $g_\ba$ on the appropriate sets.
The collection of
 $h_\ba$ for each $\ba$ gives a cochain $h$ with  
$q(h) = (h\restr{U}, h\restr{V}) = (f,g)$.
This shows exactness at the middle.

Exactness at the left end merely requires that $q$ be injective,
namely that if $f\restr{U\cap A_\ba} = g\restr{U\cap A_\ba}$
and $f\restr{V\cap A_\ba} = g\restr{V\cap A_\ba}$, then $f=g$ on 
$(U\cup V) \cap A_\ba$.  This is the second axiom of sheaves.

Exactness at the right end requries that $r$ be surjective.
This is the most interesting part, as we don't have any sheaf axioms left; 
the argument in \cite{BT} 
uses partitions of unity, which do not exist for all sheaves. 
In our case, we have partitions of unity in the $t$ directions, 
which is sufficient.

Explicitly, let $A$ be an element from the cover $\mathcal{A}$.  
An element of $\J$ looks like $h = a(t)e^{it\theta}$; suppose 
such an element is given on $A\cap U \cap V$.  
Let $\rho_U(t), \rho_V(t)$ be two functions on $\R$ 
supported on $I_U$ and $I_V$ respectively, with
$\rho_U + \rho_V \equiv 1$.
Then $\rho_V h$ is a section over $U$, and $-\rho_U h$ is a section over
$V$, both of which are in $\J$,
since they are of the form (smooth function)$\times \, e^{it\theta}$.
Furthermore, $\rho_V h - (-\rho_U h) = h$ on $U\cap V \cap A$.
Thus $\rho_V h \oplus (-\rho_U h)$ maps to $h$ under $r$.

Therefore $r$ is surjective, and the sequence \eqref{ses} is exact.
\end{proof}

\begin{prop}[Band Sheaf Mayer-Vietoris]\labell{prop:MV}
Let $U$ and $V$ be bands in $M$.  There is a long exact sequence
of sheaf cohomology
\begin{multline}\labell{eq:MV}
\cdots 
\to H^k(U\cup V, \J) 
\to H^k(U,\J) \oplus H^1(V,\J) \to \\
\to H^k(U\cap V, \J) 
\to H^{k+1} (U\cup V, \J) \to \cdots
\end{multline}
This sequence holds both for actual sheaf cohomology, 
and also for cohomology computed with respect to a particular cover.
\end{prop}

\begin{proof}
For each covering $\A$ of $M$, there is a short exact sequence of cochain 
complexes \eqref{ses}, computed with respect to $\A$.
Each  induces a long exact sequence in cohomology, as usual,
which is \eqref{eq:MV} with respect to the cover $\A$,
and so the cover-specific result is shown.

Next, the sheaf cohomology of $M$ is 
the direct limit, over the set of open covers of $M$,
of the cohomology with respect to each cover.
The maps in the directed system of open covers are just 
the restriction of sheaf elements to smaller sets in a finer cover
(see~\ref{ss:sheaf-cohom}).
Since the maps in the directed system of open covers are just 
restrictions, which commute with the coboundary maps in the 
cochain complexes, they induce transformations of cochain complexes.
By Lemma~\ref{sesnat}, this induces a map between the corresponding
long exact sequences of cohomology, namely, the Mayer-Vietoris sequences
of cohomology with respect to the particular covers.
Thus we have, in essence, a directed system of long exact sequences 
of the form \eqref{eq:MV}, \emph{with respect to the particular covers.}

The maps in this directed system, which are restrictions, 
commute with the maps in 
the long exact sequences, which are also restrictions or subtractions 
(defined following \eqref{ses}).
Hence, by Lemma~\ref{limnat}, the exactness 
passes to the direct limit, and thus we have the sequence \eqref{eq:MV}
in actual cohomology.
\end{proof}

\begin{cor}\labell{cor:MV-cov}
Let $U$ and $V$ be bands in $M$, 
such that $U \cap V$ does not contain a Bohr-Sommerfeld leaf.
Let $\B$ be a brick wall cover of $U\cup V$, which restricts to
a brick wall cover on $U\cap V$ with only two layers.  Then
\begin{equation}
H^k_\B (U\cup V, \J) \cong H^k_\B(U,\J) \oplus H^k_\B(V, \J).
\end{equation}
(The case we're most concerned with is $k=1$.  
It is true for other values of $k$ as well, though in this case 
both sides are 0.)
\end{cor}
\begin{proof}
We have the sequence \eqref{eq:MV}
with respect to the cover $\B$.
Now $H^{k-1}_\B(U\cap V, \J)$ and $H^k_\B(U\cap V, \J)$ are both zero
for all values of $k$,
by the results in Section~\ref{s:advbrickwall}:
if $k \neq 1$ then $H^k_\B(U\cap V,\J) = 0$ automatically, and
if $k=1$ then it follows from the assumption 
that $U\cap V$ contains no Bohr-Sommerfeld leaf.
Therefore, the map in \eqref{eq:MV}
\[ H^k_\B (U\cup V, \J) \to H^k_\B(U,\J) \oplus H^k_\B(V, \J) \]
is an isomorphism.
\end{proof}

\subsection{Refinements and covers: scaling the brick wall}\labell{ss:cyl-refcov}

In this section we show that the cohomology computed 
in the preceding sections 
is the actual sheaf cohomology of $U$.

\begin{thm}\labell{cohom-nocover}
Let $U$ be a band in the cylinder $M$.
Then 
\begin{equation}\labell{band-cohom}
H^k (U;\J) \cong 
\begin{cases}
  \C^\nu &k = 1\\
  0   &k\neq 1
\end{cases}
\end{equation}
where $\nu$ is the number of \BS\ leaves contained in $U$, 
and where the cohomology is the actual sheaf cohomology.
\end{thm}

The proof of this theorem uses a couple of  technical lemmas.


\begin{lemma}\labell{brick-cohom}
The cohomology 
of a band which contains no \BS\ leaf, 
with respect to a brick wall cover, is trivial, 
even if the brick wall has countably many layers.  
\end{lemma}

\begin{proof}
Suppose a brick wall covering $\B$ is given.

We first show that $H^1=0$.  This requires that, given a 1-cocycle $\eta$,
we can find a 0-cochain $\beta$ whose coboundary is $\eta$.

Index the layers by some appropriate set of integers, and 
denote the $k^\text{th}$ layer by $R^k$, made up of $m_k$ bricks.
Then $\eta$ is a collection of $m_k$ functions defined on the intersections
between bricks in each layer $R^k$, \emph{plus} a number of 
functions defined on the intersections between adjacent layers.  
The 0-cochain $\beta$, on the other hand, is a collection of $m_k$
functions $\beta^k_1,\ldots,\beta^k_{m_k}$ on each of the bricks
$B^k_1,\ldots B^k_{m_k}$ in $R^k$, for all $k$.

As we saw in section~\ref{s:advbrickwall}, however, the 
functions $\beta^k_1,\ldots,\beta^k_{m_k}$ are uniquely determined 
from the intersections within the layer $B_k$.  
We also saw that the cocycle conditions guarantee that these solutions
are consistent with the requirements that come from the parts of $\eta$ 
defined on the intersections between bricks in different layers.
More briefly, the $\beta^k_j$'s are defined entirely by the parts 
of $\eta$ living on the $k^\text{th}$ layer, and the layers 
above and below don't interfere.

Thus, given $\eta$, $\beta$ is uniquely determined on each brick,
and thus uniquely determined as a cochain.
Thus $\eta$ is a coboundary, and $H^1$ is zero.

The argument for the cohomology in other dimensions is similar.
\end{proof}

The following is a standard result (for example, 
\cite{dugundji}, Theorem XI.4.5):

\begin{lemn}[Lebesgue's Number Lemma]
Given a covering of a compact metric space $X$, 
there is a number $\delta$ such that any subspace of $X$ 
of diameter less than $\delta$ is contained in one of the sets
of the cover.
\end{lemn}

\begin{lemma}\labell{brickwall}
Any open cover $\A$ of a band $U$
has a refinement which is a brick wall cover
(with possibly countably many layers of bricks).
(Recall that one cover $\B$ is a \define{refinement} of 
another cover $\A$ if every element of $\B$ is contained 
in some element of $\A$.)
\end{lemma}

\begin{proof}
Assume a cover $\A$ of $U= I\cross S^1$ is given.  

Divide $U$ up into a countable union of \emph{closed} bands 
$\{R^k \}_{k=-\infty}^\infty$,
by dividing $I$ up into a union of closed intervals,
with the properties 
\begin{itemize}
\item each $R^k$ overlaps its two neighbours $R^{k-1}$ and $R^{k+1}$
by some positive amount (i.e.\ their interiors overlap);
\item $R^k$ does not intersect any $R^j$ other than its two neighbours; and
\item the union of the $R^k$ equals $U$.
\end{itemize}
The covering $\A$ induces a covering of each layer $R^k$.

Starting with $R^0$ and proceeding inductively in both directions, 
choose a covering $\{ B^k_1,\ldots,B^k_{n_k} \}$ of $R^k$,
with each $B^k_j$ a \emph{closed} rectangle, 
so that 
\begin{itemize}
\item each $B^k_j$ is contained in some element of $\A$
(which, if the $B^k_j$ are chosen small enough,
 is possible by Lebesgue's number lemma---if necessary, divide $R^k$ 
into more layers),
\item the collection 
$\B_k = \{ \mathring{B}^k_1,\ldots,\mathring{B}^k_{n_k} \}$
form a layer of bricks (in terms of their overlaps), 
where $\mathring{}$ denotes interior, and
\item the coverings $\B_{-k},\ldots \B_k$ form a brick wall 
(mostly, this amounts to choosing the overlaps between bricks 
in one layer to avoid the overlaps between bricks in the 
neighbouring layers).
\end{itemize}

Then the cover $\B := \bigcup_{k=-\infty}^\infty \B_k$ is a brick wall
covering of $U$ which is a refinement of $\A$.
(The fact that the partial cover is a brick wall at each stage 
ensures that the entire cover is a brick wall; 
the condition that the interiors of the $R^k$'s overlap 
ensures that there are no ``gaps'' in the cover, 
and the fact that the $R^k$'s cover $U$ guarantees that the 
$\mathring{B}^k_j$'s 
cover all of $U$.)
\end{proof}

\gap

\begin{proof}[Proof of Theorem \ref{cohom-nocover}]
Let $\B$ be a brick wall 
covering of $U$, with possibly countably many layers, 
 such that no Bohr-Sommerfeld leaf lies in more than one
layer of bricks.

Denote by $B_n$ the layer containing the $n^\text{th}$ \BS\ leaf
(i.e.\ the one corresponding to the integer $n$),
and let $V_n$ be the union of all layers of bricks between 
(but not including)
$B_n$ and $B_{n+1}$.  
Note that $B_n$ and $V_n$ are bands, and $V_n$ contains no \BS\ leaf.
Then by Lemma~\ref{brick-cohom}, $H^*_\B(V_n;\J)=0$, 
while from the calculations 
in sections~\ref{ss:sheaf-cylinder} and \ref{s:advbrickwall}, 
$H^*_\B(B_n;\J) \cong \C$ for each $n$ appearing.

Since $U$ is a finite union of $B_n$'s and $V_n$'s, $H^*_\B(U;\J)$ 
is the finite sum of the cohomologies of $B_n$ and $V_n$,
 by the cover-specific Mayer-Vietoris.  
Thus we have 
\begin{equation}
H^k_\B (U;\J) \cong 
\begin{cases}
  \C^\nu &k = 1\\
  0   &k\neq 1
\end{cases}
\end{equation}
for any brick wall covering $\B$.

However, by Lemma~\ref{brickwall}, 
any covering of $U$ has a refinement which is a brick wall.
In the language of direct limits, 
the set of brick wall coverings is cofinal in the directed system
of  coverings used to calculate sheaf cohomology.  
(See Lemma~\ref{cofinal}.)
Since the cohomology computed using any brick wall cover is the same, 
cofinality means that the actual sheaf cohomology
is the one computed with these covers, 
and thus we have \eqref{band-cohom}.
\end{proof}

\begin{cor}\labell{cor:MV}
Let $U$ and $V$ be bands in $M$, 
such that $U \cap V$ does not contain a Bohr-Sommerfeld leaf.
Then
\begin{equation}
H^k (U\cup V, \J) \cong H^k(U,\J) \oplus H^k(V, \J).
\end{equation}
in actual sheaf cohomology. 
\end{cor}
\begin{proof}
This is just Corollary~\ref{cor:MV-cov}, without the dependence on 
the cover.
Since we now know that
 $H^{k-1}(U\cap V, \J)$ and $H^k(U\cap V, \J)$ are both zero
for all values of $k$, without the cover dependence, 
the same argument as in Corollary~\ref{cor:MV-cov} shows that
\[ H^k (U\cup V, \J) \cong H^k(U,\J) \oplus H^k(V, \J). \]
\end{proof}

\section{The complex plane}\labell{C-section}

The second model space we consider is the complex plane $\C$.  
In this section we describe the set-up of the model space,
and compute the the sheaf cohomology of $\C$ by hand.
The result we obtain is a little surprising, in that the 
count over \BS\ leaves excludes the origin.  
The heart of the surprising result is Propositions~\ref{prop:zero}
and~\ref{disczero}.

The coordinates we use on $\C$ are $(s,\phi)$, 
where $(r,\phi)$ are standard polar coordinates and $s=\frac{1}{2}r^2$.
In these coordinates, the standard symplectic form has
the expression $\omega = ds\wedge d\phi$.
(Note that $d\phi$ is not defined at $s=0$, but this form extends 
smoothly there.)

The plane is equipped with a singular real polarization given by 
the distribution $P = \text{span}\{ \frac{\del}{\del \phi} \}$,
which is integrable.
Its integral manifolds are the circles of constant $s$.  
Note that this is not quite a foliation, because the manifold with $s=0$
is a point, while the others are circles, but it is a singular foliation,
as in~\ref{ss:qn-pol}.

\subsection{The sheaf of sections flat along the leaves}

Let $\LL$ be the trivial line bundle $\C \cross \C$.  
The symplectic form $\omega$ is exact
$\bigl(\omega = d(s\, d\phi)\bigr)$, and so 
the connection defined in the canonical trivialization
of $\LL$ by 
\[ \nabla_X \sigma = d\sigma(X) - \sigma i s\, d\phi(X) \]
(where $\sigma \colon \C \to \C$)
is a prequantization connection.

As before, we denote by $\J$ the sheaf of sections flat along the leaves.

\begin{prop}
The sections which are flat along the leaves are of the form 
$a(s) e^{is\phi}$, for arbitrary smooth functions $a$.
\end{prop}

\begin{proof}
The argument is very similar to the argument in section \ref{ss:BohrS},
with $t$ replaced by $s$.
A section of $\LL$ over $U$ is given by 
a map $\sigma \colon U \to \C$.  
It will be  flat along the leaves if 
\[ 0 = \nabla_X \sigma = X(\sigma) - \sigma \, i s \,d\phi(X)
	\quad \forall X \in P.  \]
Since $P$ is the span of $\frac{\del}{\del\phi}$,
this is the same as 
\[
0= \nabla_{\!\ddphi} \sigma = \frac{\del\sigma}{\del\phi}
	-\sigma i s,
\]
i.e.
\[
\frac{\del\sigma}{\del\phi} = i s \sigma.
\]
Thus, the sections flat along the leaves are those of the form
\begin{equation}\labell{Cflat}
 \sigma = a(s) e^{i s \phi}
\end{equation}
for $a$ some smooth function of $s$.
\end{proof}

This calculation applies anywhere but at the origin ($s=0$),
as $\ddth{}$ is not defined there.  
However, as we will see (Prop~\ref{prop:zero}), this is enough
to determine the value of $\sigma$ at 0.

\begin{lemma}\labell{C:BohrSset}
The Bohr-Sommerfeld leaves on $\C$ are precisely the circles 
$\{ s=k \}_{k\in\N}$.  
In other words, if $s$ is not an integer, 
there is no (nonzero) flat section defined over all of $\lf_s$.
\end{lemma}

\begin{proof}
This is exactly the same argument as in Section \ref{s:cylinder}:
in order for $e^{is\phi}$ to be defined on an entire leaf, 
that is, the entire range of $\phi$ from 0 to $2\pi$, 
it is necessary for $s$ to be an integer.  
\end{proof}

\begin{prop}\labell{prop:zero}
If $U$ is a small open disc centered at 0, then $\J(U)=0$.
\end{prop}

\begin{proof}
By Lemma \ref{C:BohrSset}, a leaf only possesses a global flat section 
if its $s$ value is an integer.  Since $U$ is a disc, it is made up 
of the union of leaves.  On every leaf which has a non-integer $s$ 
value, 
and in particular, all those near the origin,
a flat section must be zero.
Since a section is continuous,
it must be zero everywhere.  Thus there are no flat sections over
$U$ other than the zero section.
\end{proof}

Another way of stating this is: Any section 
must be 0 on any open disc around the origin
on which it is defined.
This implies that the stalk of $\J$ over 0 is zero, 
but it is actually a stronger condition.

\subsection{Cohomology}\labell{ss:C-cohom}
In this section we calculate $H^k(U,\J)$ for certain
open sets $U \subset \C$ by a similar procedure as in section~\ref{s:cylinder}.

First of all, suppose $U$ is an annulus centered at the origin, 
$\{ (s,\phi) \mid r_0 < s < r_1, \quad r_0 > 0 \}$.  
This set is analogous to the band considered in section \ref{s:cylinder},
and we have  the same result for it:

\begin{prop}\labell{prop-annulus}
If $U$ is an annulus centered at the origin which contains at most
one Bohr-Sommerfeld leaf, its sheaf cohomology is
\begin{equation}
 H^1(U, \J) = 
\begin{cases}
  \C \quad\text{if $U$ contains a Bohr-Sommerfeld leaf}\\
  0  \quad\text{if not}
\end{cases}
\end{equation}
and $H^k(U,\J) =0 $ for $k \neq 1$.
\end{prop}

\begin{proof}
This is the same proof as for Prop~\ref{e3-band-cohom}.  
The set $U$ and the elements of the sheaf have 
the same form here as there, with $t$ replaced by $s$, 
and the same argument goes through word for word.
\end{proof}

\begin{figure}[htbp]
\centerline{\includegraphics[width=1.75in,height=1.75in]{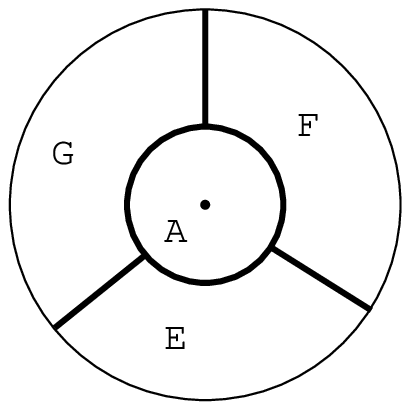}}
\caption{The cover $\A$}\label{fig1}
\end{figure}

Next, consider the open set $A$ which is a disc centred at the origin,
surrounded by three sets $E, F, G$ in a ring, none of which
intersects a Bohr-Sommerfeld orbit.  (See Figure~\ref{fig1}.)
Let $\A$ denote the cover $\{ A, E, F, G \}$, and $U =A\cup E\cup F\cup G$. 
We will calculate $H^k_\A(U;\J)$, beginning with $H^1$.

As before, we set out our bookkeeping conventions at the outset.
We will use the same conventions as in~\ref{ss:sheaf-cylinder}, namely,
that intersections will be written in the order $EF$, $FG$, $GE$. 
For intersections involving $A$, the $A$ will always come first: 
$AEF$, etc.  On an intersection, we will always use the $\phi$ 
coordinate from the set written \emph{first} in these conventions,
\emph{unless that set is $A$,} in which case we use the second set 
(since there is no $\phi$ coordinate defined on all of $A$).
We choose the $\phi$ coordinates so that 
$\phi_G = \phi_E + 2\pi$ on $GE$ while the $\phi$ 
coordinates agree on the other intersections.  
Finally, let $I_A$ denote the range of $s$-intervals covered by $A$, 
$I_O$ (for ``outer'') denote the range of $s$-intervals covered
by $E$, $F$, and $G$, and $I_{AO} = I_A \cap I_O$.

Let $\alpha = \{\alpha_{VW}\}_{V,W=A,E,F,G}$ be a 1-cochain.
The coboundary of $\alpha$ is given by 
\[ (\delta\alpha)_{VWX} = \alpha_{WX} - \alpha_{VX} + \alpha_{VW}. \]
Thus if $\alpha$ is a cocycle,  $\delta\alpha=0$ and
we have the following equations:
\begin{equation}\labell{pre-cocyc}
\begin{split}
h_{EF} e^{is\phi_E} - h_{AF} e^{is\phi_F} + h_{AE} e^{is\phi_E} =0 
\qquad \text{on } AEF\\
h_{FG} e^{is\phi_F} - h_{AG} e^{is\phi_G} + h_{AF} e^{is\phi_F} = 0
\qquad \text{on } AFG\\
h_{GE} e^{is\phi_G} - h_{AE} e^{is\phi_E} + h_{AG} e^{is\phi_G} =0 
\qquad \text{on } AGE
\end{split}
\end{equation}
where the $h$'s are smooth functions of $s$ defined on the appropriate 
interval.
Cancelling the factors of $e^{is\phi}$, we get
\begin{subequations}\labell{cocycle}
\begin{align}
h_{EF} - h_{AF} + h_{AE} &= 0 \labell{coc1}\\
h_{FG} - h_{AG} + h_{AF} &= 0 \labell{coc2}\\
h_{GE} - e^{-2\pi is} h_{AE} + h_{AG} &= 0 \labell{coc3}
\end{align}
\end{subequations}
as functions of $s$, on $I_{AO}$.

Next, to see if $\alpha$ is a coboundary, 
we ask if there exists a 1-cochain $\beta$ such that 
$\delta \beta = \alpha$.
Unravelling the definitions, 
given $h_{VW}(s)$ satisfying equations \eqref{cocycle}, we seek 
elements of $\J$ 
$\{g_Ve^{is\phi_V}\}_{V=A,E,F,G}$, each defined on the corresponding set,
such that 
\begin{equation}\labell{cobdy}
g_W e^{is\phi_W} - g_V e^{is\phi_V} = h_{VW}e^{is\phi} 
\qquad\text{on } V \cap W
\end{equation}
for all choices of $V,W$.

First, on $AE$, \eqref{cobdy} would give
\begin{equation}
g_E e^{is\phi_E} - g_A e^{is\phi_A} = h_{AE}e^{is\phi_E}
	\qquad \text{on } AE,
\end{equation}
but by Prop.\ \ref{prop:zero}, 
there are no nonzero sections on $A$.  So the term involving 
$g_A$ is zero, and we get $g_E = h_{AE}$ on $I_{AO}$.
Similarly, $g_F = h_{AF}$ and $g_G = h_{AG}$ (on $I_{AO}$).
This determines $\beta$ on all intersections with $A$.

For the definition of $\beta$ outside of $A$, we 
need to define the $g$s on the rest of $I_O$.

\begin{lemma}
The functions $g_E$, $g_F$, and $g_G$ extend to functions on all of $I_O$.
\end{lemma}

\begin{proof}
The functions $h_{EF}$, $h_{FG}$, and $h_{GE}$ are defined on $I_O$;
the functions $h_{AE}$, $h_{AF}$, and $h_{AG}$ 
are defined on $I_{AO}$.
At present, we have defined $g_E$, $g_F$, and $g_G$ only on $I_{AO}$.

Adding up the three equations 
\eqref{coc1}, \eqref{coc2}, and \eqref{coc3}
and rearranging, we obtain
\begin{equation}\labell{gfulls}
(1-e^{-2\pi is}) h_{AE}(s) = - \bigl(h_{EF}(s) + h_{FG}(s) + h_{GE}(s) \bigr)
\end{equation}
which is true on $I_{AO}$.  
But the functions on the right side are defined on all of $I_O$.
Provided $e^{-2\pi is} \neq 1$ 
(which is true by our assumption that none of these sets contain a \BS\ leaf),
this gives that $h_{AE}(s)$ on $I_{AO}$
is equal to a function that is defined on all of $I_O$.
We can use this to define  $g_E$ on the rest of $I_O$.

Similarly, by adding
$e^{-2\pi is} \cdot$\eqref{coc1}, \eqref{coc2}, and \eqref{coc3}, 
we get a similar result for $g_F$; the case for $g_G$ is similar.
For reference, we collect here the equations defining the $g$'s:
\begin{subequations}\labell{gref}
\begin{align}
(1-e^{-2\pi is}) g_E &= - \bigl(h_{EF} + h_{FG} + h_{GE} \bigr)\labell{gref1}\\
(1-e^{-2\pi is}) g_F &= - \bigl(e^{-2\pi is} h_{EF} + h_{FG} + h_{GE} 
	\bigr)\labell{gref2}\\
(1-e^{-2\pi is}) g_G &= - \bigl(e^{-2\pi is}h_{EF} + e^{-2\pi is} h_{FG} 
	+ h_{GE} \bigr)\labell{gref3}
\end{align}
\end{subequations}
\end{proof}

These extensions define $\beta$ on the rest of the sets 
$E$, $F$, and $G$.  

Finally, there are additional conditions that $\beta$ must satisfy, 
which arise from the intersections 
around the outside of the `ring,' namely $EF$, $FG$, and $GE$.  
If we apply \eqref{cobdy} to these intersections, we get the 
following conditions, which are required for $\delta \beta = \alpha$:
\begin{equation}
\begin{split}
h_{EF} e^{is\phi_E} = g_F e^{is\phi_F} - g_E e^{is\phi_E}
	\qquad \text{on } EF\\
h_{FG} e^{is\phi_F} = g_G e^{is\phi_G} - g_F e^{is\phi_F}
	\qquad \text{on } FG\\
h_{GE} e^{is\phi_G} = g_E e^{is\phi_E} - g_G e^{is\phi_G}
	\qquad \text{on } GE
\end{split}
\end{equation}
Using the convention that  $\phi_G = \phi_E + 2\pi$ and 
cancelling factors of $e^{is\phi}$ as before, we obtain
\begin{equation}\labell{ring-cobdy}
\begin{split}
h_{EF} &= g_F - g_E \\
h_{FG} &= g_G - g_F \\
h_{GE} &= e^{-2\pi is} g_E - g_G 
\end{split}
\end{equation}
on $I_{AO}$.
However, the $g$s given by \eqref{gref}, satisfy these equations
with no further restrictions (still assuming $e^{2\pi is} \neq 1$).
For example, using \eqref{gref2} and \eqref{gref1},
\begin{equation*}
\begin{split}
g_F - g_E
   &= - \tfrac{1}{(1-e^{-2\pi is})}
\bigl(e^{-2\pi is}h_{EF} + h_{FG} + h_{GE} \bigr)
   +  \tfrac{1}{1-e^{-2\pi is}}\bigl(h_{EF} + h_{GE} + h_{FG} \bigr)\\
   &=  \frac{1}{1-e^{-2\pi is}} \bigl( (1-e^{2\pi is}) h_{EF} \bigr)\\
   &= h_{EF}
\end{split}
\end{equation*}
as required.

(In fact, the required conditions \eqref{ring-cobdy} 
are the same as the cocycle conditions \eqref{cocycle}, 
except that the latter only apply on 
the intersections with $A$, while the above calculation applies on the 
entire domain of $\beta$.  
This is reminiscent of what we saw in \ref{s:advbrickwall}
(around equations~\eqref{cij} and \eqref{cij2}),
where the conditions for compatibility of two sets of solutions
were precisely the cocycle conditions.)

Therefore, given any cocycle $\alpha \in C^1$, there exists a 
$\beta \in C^0$ such that $\delta \beta = \alpha$, and hence
every 1-cocycle is a coboundary.  
Thus, we have proved the following:
\begin{thm}\labell{disczero}
For $U \in \C$ an open disc centred at 0, and $\A$ the cover given
in Figure~\ref{fig1}, 
 $H^1_\A(U;\J) = 0$.
\end{thm}

The cohomology in all other dimensions is trivial: 

\begin{prop}\labell{p:C-cohom-k}
 $H^k_\A (U;\J) = 0$ for $k\neq 1$ as well.
\end{prop}

\begin{proof}
First consider $k=0$.  
As noted in section~\ref{ss:sheaf-cylinder}, $H^0(U;\J)$ is the space of global
sections of $\J$ over $U$.  
Since there are no such sections, $H^0$ is zero.

Next, consider $H^2$.

A 2-cochain $\alpha$ is a collection 
$ \{ a_{AEF} e^{is\phi_E}, a_{AFG} e^{is\phi_F}, 
a_{AGE} e^{is\phi_G} \}$,
and is automatically a cocycle.

A 1-cochain $\beta$
is a collection $\{ b_{VW} e^{is\phi_U} \}_{U,V = E,F,G,A}$.
Now $\alpha$ will be $\delta \beta$ if 
\begin{equation}
\begin{split}
b_{EF} e^{is\phi_E} - b_{AF} e^{is\phi_F} + b_{AE} e^{is\phi_E} 
	= a_{AEF} e^{is\phi_E} \quad \text {on } AEF \\
b_{FG} e^{is\phi_F} - b_{AG} e^{is\phi_G} + b_{AF} e^{is\phi_F} 
	= a_{AFG} e^{is\phi_F} \quad \text {on } AFG \\
b_{GE} e^{is\phi_G} - b_{AE} e^{is\phi_E} + b_{AG} e^{is\phi_G}
	= a_{AGE} e^{is\phi_G} \quad \text {on } AGE 
\end{split}
\end{equation}
Using the convention that $\phi_G = \phi_E + 2\pi$ on $GE$, 
and cancelling factors of $e^{is\phi}$, we obtain:
\begin{equation}\labell{C-2-coc-system}
\begin{split}
a_{AEF} &= b_{EF} - b_{AF} + b_{AE} \\
a_{AFG} &= b_{FG} - b_{AG} + b_{AF} \\
a_{AGE} &= b_{GE} - e^{-2\pi is} b_{AE} + b_{AG} 
\end{split}
\end{equation}
on $I_{AO}$.
This is a system of 3 equations in 6 unknowns, and thus has many solutions.
(It is easy to check that it is consistent.)  
There is the question of extendability: for example, $a_{AGE}$ is only
defined for $s \in I_{AO}$, while $b_{EF}$ needs to be 
defined for $s \in I_O$.  
However, if we set $b_{EF}$, $b_{FG}$, and $b_{GE}$ to zero, 
the resulting system in $b_{AE}$, $b_{AF}$, and $b_{AG}$ 
is the same system as \eqref{coc-system}, which has a solution
provided that $e^{2\pi is} \neq 1$.  Since there are no Bohr-Sommerfeld
leaves in any of the sets we are considering, this condition holds, 
and so the sytem \eqref{C-2-coc-system} always has a solution.  
Thus, any cocycle is a coboundary, and the second cohomology is trivial.

Finally, for $k \geq 3$, there are no $(k+1)$-fold intersections, 
and so there is no cohomology.
\end{proof}

\begin{prop}
The cohomology calculated with respect to this particular cover 
is the actual sheaf cohomology.
\end{prop}

\begin{proof}
This is essentially the same argument as in Section~\ref{ss:cyl-refcov}.
The only modification that we need is in the form of the cover, 
as a brick wall doesn't immediately apply to the complex plane.
The types of covers we use are as follows: a brick wall covering 
(a brick wall in $s$, $\phi$ coordinates) of 
the set $\{ s\neq 0\}$, plus a disc centred at zero.  
It is clear that any cover has a refinement of this form: 
away from the origin the same argument applies as with the brick wall, 
and around the origin we need only take a small enough disc.
\end{proof}

\subsection{Mayer-Vietoris}\labell{ss:C-MV}

\begin{prop}
If each of $U,V\subset \C$ is either an open annulus or an open disc, 
centered at the origin, then the sequence 
\begin{multline}
\cdots 
\to H^1(U\cup V, \J) 
\to H^1(U,\J) \oplus H^1(V,\J) \to \\
\to H^1(U\cap V, \J) 
\to H^2 (U\cup V, \J) \to \cdots
\end{multline}
as in Prop.\ \ref{prop:MV}, is exact.
\end{prop}

\begin{proof}
Away from 0, the exact same argument as in section \ref{ss:cyl-MV}
applies, with $t$ replaced by $s$.
Near 0, we need another argument.

Assume $U$ is a disc centered at 0, and $V$ is an annulus 
centered at 0 which overlaps with $U$ and whose closure 
does not contain 0.
Consider the sequence of \v Cech cochain complexes
computed with respect to some cover $\B$
(which here may be any cover):
\begin{equation}\labell{C:ses}
0 \to C^*_\B (U \cup V, \J) \stackrel{q}{\to}
C^*_\B (U,\J) \oplus C^*_\B (V,\J) \stackrel{r}{\to}
C^*_\B (U \cap V, \J) \to 0
\end{equation}
which is \eqref{ses} from section~\ref{ss:cyl-MV}
(where the definition of the maps $r$ and $q$ is given).
Exactness at the left and the centre only relies on properties of 
sheaves, as described in the proof of Prop~\ref{p:exact},
and thus still hold in this case.  
Exactness on the right requires $r$ to be surjective, 
namely, given some $h \in \J(B\cap U\cap V)$, we need 
$f \in \J(U)$ and $g \in \J(V)$ whose difference on 
$B\cap U\cap V$ is $h$.
As in the proof of Prop~\ref{p:exact}, we take a partition
of unity $(\rho_U(s),\rho_V(s))$ over the $s$-intervals covered 
by $U$ and $V$, and let $f=\rho_Vh$ on $U\cap B$, $g=\rho_Uh$ on $V\cap B$.
The only possible concern is that, if $B$ contains the origin,
$f$ must be zero on any disc centred at 0 contained in $B\cap U$.
Since $\overline{V}$ does not contain 0, this means there 
is a disc around zero outside $\overline{V}$;
and since $\rho_V$ is zero outside $\overline{V}$,
$f$ must therefore be zero on this disc, and therefore 
on any disc in $B\cap U$.
\end{proof}

By repeated applications of Mayer-Vietoris, we obtain:
\begin{thm}
Let $U\subset \C$ be either an annulus or a disc, centred at the origin.
Then 
\[ H^1(U,\J) \cong  \C^m \]
where $m$ is the number of \BS\ leaves, \emph{excluding} 
the origin, contained in $U$.
\end{thm}

\section{Example: $S^2$}\labell{s:s2}

\begin{figure}[h]
\centerline{\includegraphics{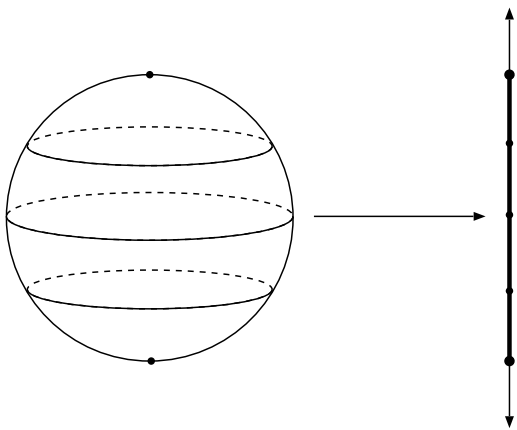}}
\caption{The moment map on the 2-sphere}\label{fig:ubiquitous}
\end{figure}

The simplest example of a toric manifold, which is ubiquitous\footnote{
This is probably because it is the only toric manifold of dimension less 
than 4 (by Delzant's classification, see~\cite{ACL}), and so it's 
the only one that can be drawn on a page.}
in textbooks on symplectic geometry, is $S^2$ with an action of 
$S^1$ by rotations about the $z$-axis.  
In this case the moment map is simply the height function.
(See Figure~\ref{fig:ubiquitous}.)
This example serves as a good illustration of the results of this paper:
despite its simplicity, it contains the essential idea of our method.

\begin{figure}[htbp]
\end{figure}

\begin{figure}[h]
\centerline{\includegraphics{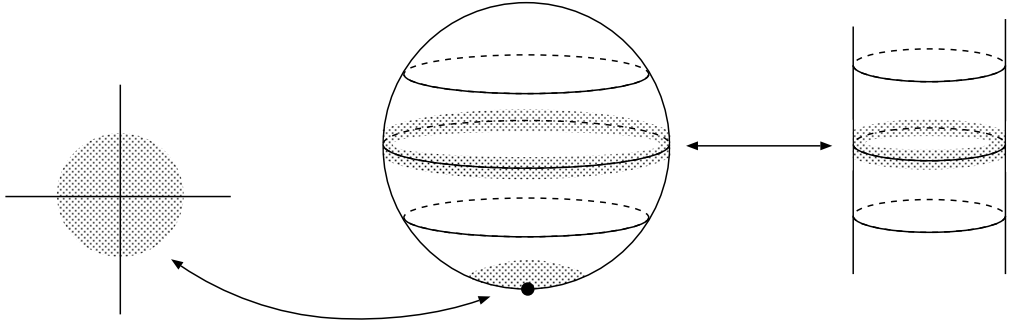}}
\caption{Neighbourhoods of orbits on the 2-sphere}\label{fig:s2}
\end{figure}

In fact, once we draw the picture in Figure~\ref{fig:s2},
further explanation is almost unnecessary.
The orbits of the circle action are circles of constant height, plus the 
two singular orbits at the north and south poles.  
A neighbourhood of a circle orbit looks like a neighbourhood of a circle 
in the cylinder, and a neighborhood of one of the poles looks like
a neighbourhood of the origin in $\R^2$, as illustrated in Figure~\ref{fig:s2}.
We have determined the cohomology of each of these neighbourhoods in 
the preceding two sections, and we transfer these results over to the sphere.

\begin{figure}[htbp]
\end{figure}

More formally, 
let $M=S^2$, and $\omega_1$ be the standard symplectic (area) form,
normalized so that the total area of the sphere is $2\pi$.
Let $\LL_1$ be the complex line bundle with Chern class 1, which is a
prequantum line bundle for $S^2$.
Let $\LL_k = \LL_1^{\tensor k}$, which will then be a prequantum line 
bundle for $S^2$ with symplectic form $\w_k := k\, \w_1$, 
by the additivity of the Chern class.

As mentioned above, the height function $\mu$ is a moment map for the 
circle action, normalized so the total height of the sphere is 1; 
similarly, if we take the same circle action but use 
the symplectic form $\omega_k$, then its moment map is 
$k$ times the height function, which we denote by $\mu_k \colon S^2 \to \R$.
The fibres give a singular Lagrangian fibration.

The Bohr-Sommerfeld leaves in this example are the 
leaves on which the moment map has an integer value,
which we can see as follows.
As in section~\ref{s:hol}, \BS\ fibres are those with trivial holonomy,
i.e.\ those leaves $\lf$ such that 
\[ \exp \left( i \int_\lf \Theta_k \right) = 1 \]
where $\Theta_k$ is a potential 1-form for the connection on the 
prequantum line bundle $\LL_k$.  
This will be true if{f} 
\[ \int_\lf \Theta_k = \int_\Sigma \omega_k \in 2\pi\Z \]
where $\Sigma$ is a surface whose boundary is $\lf$.  
Taking $\Sigma$ to be the ``bottom cap,'' the set of points 
with height less than $\lf$, then $\lf$ will be \BS\ if{f}
$k$ times the area of $\Sigma$ is
an integer multiple of $2\pi$ --- that is,
if{f} $\mu_k(\lf)$ is an integer.
(See Figure~\ref{fig:ubiquitous}, where $k=4$.
It is true in general that the \BS\ points are integer points of the moment 
map --- see Proposition~\ref{intlattice}.)
Since the image of $\mu_k$ will be the interval $[0,k]$, 
there will be $k+1$ integer points in this interval, 
including the endpoints.

Cover $\R$ by a sequence of pairwise overlapping intervals.  
This  induces a covering of $S^2$ by sets $U_j$ 
which are the equivalent in $S^2$ of bands in the model spaces.
In fact, we have a symplectomorphism from such a set to a band in 
the cylinder, or a disc in $\C$.  
By the results in Section~\ref{s:piecing}, the cohomology of $U_j$ 
is isomorphic to the cohomology of the corresponding band 
or disc in the model space.
This cohomology is trivial in all dimensions other than 1, 
as we have seen, and 
even in dimension 1 is only non-trivial if $U_j$ contains 
a non-singular \BS\ leaf.  
Adding up the results from all of the $U_j$ by Mayer-Vietoris 
(Proposition~\ref{mult-MV} below), we have:

\begin{thm}
The sheaf cohomology $H^q(S^2,\J)$ is zero if $q \neq 1$, and has 
dimension equal to the number of non-singular Bohr-Sommerfeld
leaves if $q=1$.
\end{thm}

Thus, when the prequantum line bundle is $\LL_k$, the 
quantization has dimension $(k-1)$.

\section{The multidimensional case}\labell{s:multdm}
For the case of higher dimensions, the model space we use is $\msp$, 
where $k$ will be determined by the dimension of the leaf.
In this section we describe the basic set-up of the model space, 
as well as prove some results about piecing together in the 
multidimensional case.  Computation of the cohomology of the model 
space is postponed until the following section.

\subsection{The model space}
Let $M_0 = (\R\cross S^1)^m \cross \C^k$ 
(where for the moment we write $m$ for $n-k$, simply for ease of notation),
with coordinates 
\[ (t_1,\theta_1,\ldots,t_m,\theta_m,s_1,\phi_1,\ldots,s_k,\phi_k) \]
using the same conventions as in sections \ref{s:cylinder} and \ref{C-section}.
In these coordinates, the standard symplectic form is given by
\[ \omega = dt_1 \wedge d\theta_1 + \cdots + dt_m \wedge d\theta_m 
	+ ds_1 \wedge d\phi_1 + \cdots ds_k \wedge d\phi_k, \]
which is equal to 
\begin{equation}\labell{Theta}
 d( t_1 \,d\theta_1 + \cdots + s_k \, d\phi_k). 
\end{equation}
The calculations, at least in the beginning, are exactly
the same as in the one-dimensional case, except with more indices.

The polarization is 
$P=\text{span}\{\ddth{i},\ddph{j}\} $
whose leaves are surfaces of constant $t$ and $s$.
If we map $M_0$ to $\R^{m+k}$ by projecting onto the $t$ and $s$ coordinates,
the polarization is given by the level sets of this map.
Let $\LL_0$ be the trivial bundle $M_0 \cross \C$, which we
make a prequantum line bundle by giving it a connection
whose potential 1-form is 
$\Theta =  t_1 \,d\theta_1 + \cdots + s_k \, d\phi_k$,
so that $d\Theta = \w$.

\begin{defnn}\labell{def:band}
A \define{band} in $(\R\cross S^1)^m$ is a set of the form $I\cross (S^1)^m$ 
where $I$ is an open rectangle, i.e.\ the product of intervals, 
in $\R^m$.
More generally, a \define{band}
in $(\R\cross S^1)^{n-k} \cross \C^k$
is a band in the preceding sense 
in $(\R\cross S^1)^{n-k}$,
times the product of discs centred at 0 in $\C^k$.
Even more generally, a \define{band} in a manifold $M$ is a set 
symplectomorphic to a band (in the preceding sense)
in $(\R\cross S^1)^{n-k} \cross \C^k$.
\end{defnn}

\subsection{The flat sections}

\begin{prop}
A section of $\LL_0$ which is flat along the leaves locally has the form 
\begin{equation}
\sigma = a(t_1,t_2,\ldots,s_k)
e^{i(t_1\theta_1 + t_2\theta_2 + \cdots + s_k\phi_k)},
\end{equation}
which we may write as 
\[ a(\bt,\bs)\, e^{i\ip{\bt}{\bth}} e^{i\ip{\bs}{\bph}} \]
where $a$ is a smooth function.
\end{prop}

\begin{proof}
Let $\sigma \colon U \to \C$ be a section of $\LL_0$, 
which we view as a $\C$-valued function using the
 canonical trivialization.
Then $\sigma$ is flat along the leaves if 
$\nabla_X \sigma = 0$ for all $X\in P$; this will be true iff
\begin{equation}
\begin{split}
\nabla_{\ddth{j}}\sigma = 0 \\
\nabla_{\ddph{l}}\sigma = 0
\end{split}
\end{equation}
for all appropriate values of $j$ and $l$, 
since the connection is linear in $X$.

Applying the argument from the proof of Proposition~\ref{cyl-flat} 
to each $\theta$ and $\phi$ component in turn 
(note that the connection potential 1-form has the same form in 
each component here as it did in the proposition), we 
obtain the differential equations 
\begin{equation}
0 = \frac{\del\sigma}{\del\theta_j} - \sigma i t_j
\qquad 
0 = \frac{\del\sigma}{\del\phi_j} - \sigma i s_j
\end{equation}
for all $\theta$ and $\phi$ coordinates.
Together, these equations imply the desired result.
\end{proof}

As before, we have 
\begin{prop}
The Bohr-Sommerfeld set of $M_0$ is $\Z^m \cross \N^k \subset \R^{m+k}$.
\end{prop}

\begin{proof}
The \BS\ points are those with integer $t$ values, and positive 
integral $s$ values, by exactly the same argument as in the cases of the 
cylinder and the complex plane.
\end{proof}

\subsection{Multidimensional Mayer-Vietoris}\labell{s:mult-MV}

\begin{prop}\labell{mult-MV}
Let $U$ and $V$ be subsets of $M$ which are each the union of 
leaves of the polarization.  
Then the sequence 
\begin{multline}\labell{eq:mult-MV}
\cdots 
\to H^1(U\cup V, \J) 
\to H^1(U,\J) \oplus H^1(V,\J) \\
\to H^1(U\cap V, \J) 
\to H^2 (U\cup V, \J) \to \cdots
\end{multline}
of sheaf cohomology is exact.
\end{prop}

\begin{proof}
Let $\A$ be a cover of $U\cup V$, which induces a cover 
of $U$, $V$, and $U\cap V$ as in section~\ref{ss:cyl-MV}.
Consider the sequence 
\begin{equation}\labell{mult-ses}
0 \to C^*_\A (U \cup V, \J) \stackrel{q}{\to}
C^*_\A (U,\J) \oplus C^*_\A (V,\J) \stackrel{r}{\to}
C^*_\A (U \cap V, \J) \to 0
\end{equation}
of \v Cech cochain complexes computed with respect to the 
cover $\A$.  We claim this sequence is exact.

The argument to show this is very similar to the proof of 
Proposition~\ref{p:exact}, which is the same result in the case
of the cylinder.  
Exactness at the left and the middle follows from properties of 
sheaves, exactly as in that proof.

Exactness at the right is less straightforward.  
To show exactness, we need surjectivity of $r$; 
so let $A$ be a set in the cover $\A$ and $h$ an element of 
$\J(A\cap U\cap V)$.  We require a flat section $f$ on $U\cap A$ 
and $g$ on $V\cap A$ whose difference on $U\cap V\cap A$ is $h$.

Since $U$ and $V$ are unions of leaves, they can be 
written as $\pi^{-1}(B_U)$ and $\pi^{-1}(B_V)$ for some 
subsets $B_U$ and $B_V$ of $B$.  Since $B$ is a manifold with corners,
we can find a partition of unity on $B_U$, $B_V$.  
By composing with the projection $\pi \colon M \to B$, 
we obtain a partition of unity $\rho_U$, $\rho_V$ for $U$ and $V$.  
(We take supp$(\rho_U) \subset U$.)

Now let $f = \rho_V h$ and $g = -\rho_U h$.  Then $f$ is a section
over $U$, which is flat along the leaves because $h$ is flat along 
the leaves and $\rho_V$ is constant on the leaves.  
Similarly, $g$ is a leafwise flat section defined over $V$.
It is clear that $f-g = h$ on $U\cap V\cap A$.

Therefore $r$ is surjective, and the sequence \eqref{mult-ses} is exact.

The short exact sequence \eqref{mult-ses} gives a long exact sequence
in cohomology, as usual, which is the sequence \eqref{mult-MV}
with respect to the cover $\A$.
Doing this for each $\A$ yields such a sequence for every cover.
By the same argument as in \eqref{prop:MV}, 
the exactness passes to the limit, and we have \eqref{mult-MV}.

\end{proof}

\begin{rmk}
The existence of the exact sequence~\eqref{mult-MV}
seems to be true for cohomology using arbitary sheaves---see 
\cite{Iv}, III.5.10.  
However, it is instructive to see how the particular properties 
of our sheaf allow a direct proof of Mayer-Vietoris.
\end{rmk}

\section{A better way to calculate cohomology}\labell{s:leray}

\epigraph{He that goeth forth and weepeth, bearing precious seed, shall
doubtless come again with rejoicing, bringing his sheaves with him.}
{Psalm 126:6}

For the higher dimensional model spaces, the challenge with trying to 
compute the cohomology of a cover directly, as we did in 
sections~\ref{s:cylinder} and \ref{C-section}, 
is that in order to adequately cover a set of higher dimensions, 
the covers, and thus the bookkeeping, become more and more complicated.  
Instead, 
we use a sheaf theoretic argument to obtain the cohomology.
I am grateful to Ruxandra Moraru for suggesting this
approach, and for explaining much of the sheaf theory to me.

The structure of this section is as follows.  
First, in~\ref{ss:theory}, we discuss the theoretical tools we will use.
Then, in section~\ref{sss:leray1d}, we apply them to the 2-dimensional
case (which we have already calculated), for practice.
After outlining the upcoming calculations in~\ref{leray-outline}, 
in section~\ref{leray-nd} and \ref{leray-ndsing}
we apply our tools to the higher-dimensional case.

\subsection{Theory}\labell{ss:theory}

\subsubsection{Spectral sequences}\labell{sss:ss}

We will not attempt to describe the theory of spectral sequences
in detail here, 
but refer the reader to \cite{BT}.
We will briefly review some of the facts about spectral sequences 
that will be necessary for our calculations, which we will not 
attempt to state in full generality, but just enough to suffice
for our purposes.

\begin{figure}[htbp]
\centerline{\includegraphics[width=3.5in,height=2.0in]{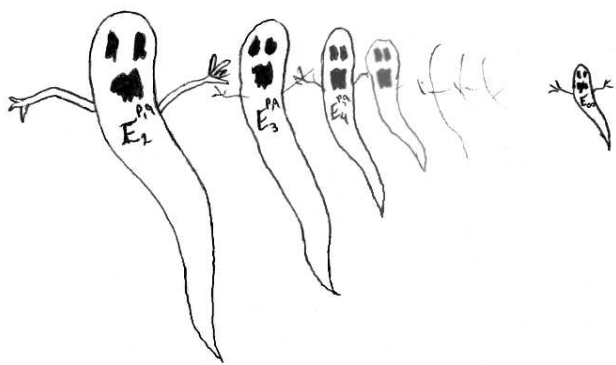}}
\caption{A spectral sequence}\label{fig:spectral}
\end{figure}

Recall that a \define{spectral sequence} is a collection
$\{ (E_r, d_r) \}$, where the $E_r$ are vector spaces,
the $d_r \colon E_r \to E_r$ are differentials (i.e.\ $d_r \circ d_r = 0$),
and each group is the cohomology of the previous one, with respect
to its differential: $E_{r+1} = H(E_r, d_r)$.
Usually, $E_r$ comes with a bigrading, and 
$d_r$ shifts the bidegree, mapping $E_r^{p,q}$ to $E_r^{p+r,q-r+1}$.

A spectral sequence is often drawn in a chart, as in Figure~\ref{fig:ss},
where each group $E^{p,q}$ is put in the appropriate square,
and where we can think of each 
different $r$ sitting on a different ``page'' in the diagram.
The maps $d_r$ go between the groups as shown; as $r$ increases,
the target of $d_r$ for a fixed source moves down the diagonal.

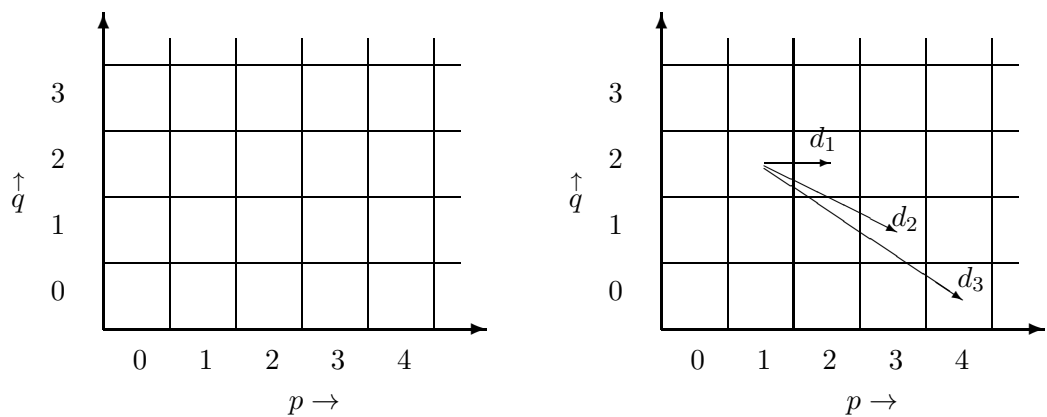
\begin{figure}[h]
\begin{center}
\makebox[150pt][r]{
\begin{picture}(200,170)(-40,-10)
\thicklines
\put(0,0){\vector(1,0){145}} 
\put(0,0){\vector(0,1){120}} 
\thinlines
\multiput(25,0)(25,0){5}{\line(0,1){110}}
\multiput(0,25)(0,25){4}{\line(1,0){135}}
\put(11,-15){0}
\put(36,-15){1}
\put(61,-15){2}
\put(86,-15){3}
\put(111,-15){4}
\put(70,-30){$p \to$}  
\put(-20,11){0}
\put(-20,36){1}
\put(-20,61){2}
\put(-20,86){3}
\put(-35,45){$\stackrel{\uparrow}{q}$} 
\end{picture}
} 
\makebox[200pt][l]{
\begin{picture}(200,170)(-40,-10)
\thicklines
\put(0,0){\vector(1,0){145}} 
\put(0,0){\vector(0,1){120}} 
\thinlines
\multiput(25,0)(25,0){5}{\line(0,1){110}}
\multiput(0,25)(0,25){4}{\line(1,0){135}}
\put(11,-15){0}
\put(36,-15){1}
\put(61,-15){2}
\put(86,-15){3}
\put(111,-15){4}
\put(70,-30){$p \to$}  
\put(-20,11){0}
\put(-20,36){1}
\put(-20,61){2}
\put(-20,86){3}
\put(-35,45){$\stackrel{\uparrow}{q}$} 
\put(39,63){\vector(1,0){25}}
\put(56,69){$d_1$}
\put(39,62){\vector(2,-1){50}}
\put(87,39){$d_2$}
\put(39,61){\vector(3,-2){75}}
\put(112,16){$d_3$}
\end{picture}
}
\end{center}
\caption{The diagram of a spectral sequence, 
and its differentials}\labell{fig:ss}
\end{figure}

If, for some $r$, all of the differentials are zero, then 
$E_{r+1} = E_r$, since $\ker d_r$ is everything and im$\,d_r$ is zero.
If there is some $s$ such that all $d_r$ are zero for $r>s$,
then all of the $E_r$'s for $r>s$ are the same, and we say the 
spectral sequence \define{stabilizes} or \define{converges}.
We denote the common value of $E_r$ for $r>s$ by $E_\infty$.

In many applications, a bigraded spectral sequence is used to obtain
some singly-graded object.  
The grading on a spectral sequence is obtained from the bigrading 
by summing along the diagonal:
\[ E^k_r = \bigoplus_{p+q=k} E_r^{p,q} \]
If we say a spectral sequence converges to some singly-graded object,
we mean it in this sense.

We observe that if, for some $r\geq 2$, the diagram of 
a spectral sequence has only one non-zero row,
then the spectral sequence stabilizes at that value of $r$.
More precisely:

\begin{prop}\labell{ss-triv}
Let $(E_r^{p,q},d_r)$ be a spectral sequence.  
Suppose that there is some number $m$ such that, 
for some $s \geq 2$, 
\[ E_s^{p,q} = 0 \qquad \text{for all } q\neq m. \]
Then the spectral sequence stabilizes for $r=s$,
i.e. $E_\infty = E_s$.
\end{prop}

\begin{proof}
Since $d_s$ maps from $E_s^{p,q}$ to $E_s^{p+s,q-s+1}$, 
if $s\geq 2$, $d_s$ maps to a group with a different value of $q$.
If there is only one value of $q$ for which the groups are non-zero,
this implies that all the differentials must be zero.
(See Figure~\ref{ss:nonzero}, where all the blank boxes are zero,
and a $\ast$ represents something possibly nonzero.)
Thus the spectral sequence stabilizes.
\end{proof}

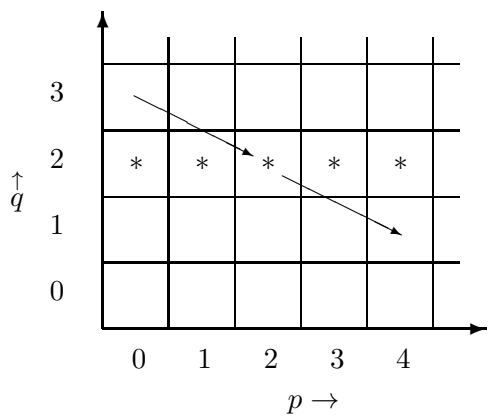
\begin{figure}[ht]
\begin{center}
\begin{picture}(200,130)(-40,-10)
\thicklines
\put(0,0){\vector(1,0){145}} 
\put(0,0){\vector(0,1){120}} 
\thinlines
\multiput(25,0)(25,0){5}{\line(0,1){110}}
\multiput(0,25)(0,25){4}{\line(1,0){135}}
\put(11,-15){0}
\put(36,-15){1}
\put(61,-15){2}
\put(86,-15){3}
\put(111,-15){4}
\put(70,-30){$p \to$}  
\put(-20,11){0}
\put(-20,36){1}
\put(-20,61){2}
\put(-20,86){3}
\put(-35,45){$\stackrel{\uparrow}{q}$} 
\put(10,60){$\ast$}
\put(35,60){$\ast$}
\put(60,60){$\ast$}
\put(85,60){$\ast$}
\put(110,60){$\ast$}
\put(68,58){\vector(2,-1){45}}
\put(12,88){\vector(2,-1){45}}
\end{picture}
\end{center}
\caption{The differentials are all zero.}\labell{ss:nonzero}
\end{figure}

\subsubsection{Leray spectral sequence}

The spectral sequence we will be using is the following.

\begin{thm}[\cite{Gode}, Theorem 4.17.1]\labell{leray}
Let $f \colon M \to B$ be a continuous map, and 
let $\sh$ be a sheaf over $M$.
Then there is a spectral sequence 
(called the \define{Leray spectral sequence}),
whose $E_2$ term is given by
\begin{equation}
E_2^{p,q} = H^p (B; R^q f_* \sh) 
\end{equation}
(where $R^q f_* \sh$ are the direct image sheaves defined below)
and which converges to $H^{p+q}(M,\sh)$.
\end{thm}

\begin{defnn}[\cite{GH}, p.\ 463]\labell{Rqsheaf}
Let $f \colon X \to Y$, and let $\sh$ be a sheaf on $X$. 
The \define{$q^\text{th}$ direct image sheaf} $R^qf_\ast\sh$
is the sheaf arising from the presheaf
\begin{equation}
 U \mapsto H^q \bigl( f^{-1}(U), \sh \bigr).
\end{equation}
In the case $q=0$, the sheaf is denoted simply by $f_*\sh$ and 
called the \define{pushforward sheaf} of $\sh$ by $f$.
\end{defnn}

In the cases we consider, rather than find the higher direct image 
sheaves directly 
by using the above definition, we will find the \emph{stalks} 
of the sheaves.  
The following result will be very useful in these calculations.

\begin{thm}[``Sheaf Theorist's Trick'']\labell{stt}
Let $f\colon X \to Y$ be a proper map between locally compact spaces,
and let $\sh$ be a sheaf on $X$.  
For $y \in Y$ and $q \in \N$ the restriction map
\begin{equation}
(R^q f_* \sh)_y \to H^q \bigl( f^{-1}(y);\sh\restr{f^{-1}(y)} \bigr) 
\end{equation}
(where $(R^q f_*\sh)_y$ denotes the stalk of the sheaf $R^q f_*\sh$ at $y$,
as in~Definition~\ref{d:stalk})
is an isomorphism.
\end{thm}

\begin{proof}
This is Theorem III.6.2 from \cite{Iv}.  See also \cite{Gode}, Remarque 4.17.1.
\end{proof}

In the cases we consider, the differentials will be trivial, and so 
usually
\begin{equation}\labell{principle2}
H^{m}(M,\sh) = \bigoplus_{p+q=m} H^p(B,R^q \pi_* \sh). 
\end{equation}
We will use the Leray spectral sequence \emph{twice} in the course of 
this calculation:
once with the map 
$(\R\cross S^1)^n \to \R^n$, and once with $(S^1)^n \to (S^1)^{n-1}$.

\begin{notn}
We generally use $\J$ to denote the sheaf of sections flat along the 
leaves.  
In what follows, use $\J_n$ to denote this sheaf over $(\R\cross S^1)^n$,
when we need to be specific about the dimension.
We will often need to consider one component of $(\R\cross S^1)^n$
at a time, with other components fixed.  
We will use the notation $\bt$ to mean $(t_1,\ldots,t_n)$ as usual,
but we will also use $\bt^{n-1}$ to mean $(t_1,\ldots,t_{n-1})$, 
to emphasize that we are not dealing with the coordinate $t_n$.  
We will use a similar notation $\bth^{n-1}$.
\end{notn}

\subsubsection{Skyscraper sheaves}

\begin{defn}  A \define{skyscraper sheaf} supported at a point $p$ 
is a sheaf $\sh$ whose every stalk is zero except the stalk at $p$.
More precisely, let $A$ be some abelian group.  Then $\sh(U) = A$ 
if $p \in U$, and $\sh(U) = 0$ otherwise.  The group $A$ is called 
the \define{tower} of the sheaf.

More generally, we allow a skyscraper sheaf to be supported at 
more than one point, provided the set of such points is discrete.  
Thus, for us, a skyscraper sheaf is one whose every stalk is zero outside of 
some discrete set.
\end{defn}

\begin{lemma}[Cohomology of a skyscraper sheaf]\labell{l:skys}
Suppose $\sh$ is a skyscraper sheaf supported on the discrete set 
$I \subset X$, with towers $A_i$ respectively.  Then the cohomology
of $\sh$ is 
\begin{equation}
\begin{split}
H^0(X,\sh) &= \bigoplus_{i\in I} A_i\\
H^q(X,\sh) &= 0 \qquad q > 0
\end{split}
\end{equation}
\end{lemma}

\begin{proof}
This is a standard result.
Since $H^0$ is just the global sections of the sheaf, the result 
for $q=0$ is immediate from the definition; the result for $q>0$ is, 
for example, Proposition IX.4.3 in~\cite{miranda}.
\end{proof}

\subsection{The case of one dimension}\labell{sss:leray1d}

Consider the map $R\cross S^1 \stackrel{\pi}{\to} \R$. 
Recall that elements of the sheaf $\J$ locally have the form
 $a(t) e^{it\theta}$.

Let $I\subset \R$ be an open interval.
By Definition~\ref{Rqsheaf}, the pushforward sheaf $\pi_*\J$ is
given by 
\[\pi_*\J(I) \cong H^0 (I\cross S^1,\J). \]
This is just the set of global sections of $\J$ over $ I \cross S^1$;
as we argued in Section~\ref{ss:sheaf-cylinder}, 
this is 0 for all intervals $I$.

Similarly, the higher direct image sheaves are given by 
\[ (R^q \pi_* \J)(I) = H^q(I \cross S^1, \J). \]
The cohomology of $I\cross S^1$ 
was computed in Section~\ref{ss:sheaf-cylinder}.
From Proposition~\ref{band-cohom}, we have that
it is zero for all $q$, except when the interval 
$I$ contains an integer, in which case $H^1 \cong \C^m$, 
where $m$ is the number of integers in $I$  
(and all other $H^q$ are still zero).
Thus:

\begin{lemma}\labell{RqpiJ1}
The $qth$ direct image sheaf $R^q\pi_*\J_1$ is 0 in the case $q\neq 1$, 
and a skyscraper sheaf supported on $\Z$, with all towers $\C$, if $q=1$.
\end{lemma}

\begin{proof}
The case $q\neq 1$ is a restatement of the preceding paragraph.  
From the same discussion, supposing $I$ is small,
\begin{equation}
R^1\pi_*\J_1(I) = 
\begin{cases}
\C &\quad \text{if } I \text{ contains an integer}\\
0  &\quad \text{if } I \text{ contains no integer}
\end{cases}
\end{equation}
This is just the definition of a skyscraper sheaf.
Its support is the set of integers.
\end{proof}

\subsection{The structure of the coming calculation}\labell{leray-outline}
In the following section, we will calculate $H^m(U;\J)$ 
for $U \subset (\R\cross S^1)^n$.  
The calculation is somewhat complicated, and so we outline it here.
\begin{enumerate}
\item\labell{one}
Let $\pi \colon (\R\cross S^1)^n \to \R^n$ be the obvious projection,
and let $\pi(U)=V$.  The Leray spectral sequence gives 
$H^m(U;\J)$ in terms of $H^p(V;R^q\pi_*\J)$, 
so we need $R^q\pi_*\J$.  

\item\labell{two}
The sheaf theorist's trick gives a stalk of $R^q\pi_*\J$ as
$H^q(\pi^{-1}(\bt);\J')$, where $\bt \in R^n$.  
Since $\pi^{-1}(\bt) = \{ \bt \} \cross T^n$,
this amounts to finding $H^q(T^n;\J')$.

\item\labell{three} 
Let $\rho\colon T^n \to T^{n-1}$ be projection onto the first $(n-1)$ 
coordinates.  The Leray spectral sequence applied to $\rho$ gives us
$H^m(T^n;\J')$ in terms of $H^p(T^{n-1};R^q\rho_*\J')$.

\item\labell{four} 
The sheaf theorist's trick gives a stalk of $R^q\rho_*\J'$
as $H^q(\rho^{-1}(x);\J')$, where $x \in T^{n-1}$.

\item\labell{five} 
Finally, $\rho^{-1}(x) = \{ x \} \cross S^1$, and so the 
calculation of $H^q(\rho^{-1}(x))$ is reduced to 
the calculation of $H^q(S^1)$, which is essentially the same calculation
as the one in section~\ref{ss:sheaf-cylinder}.

\end{enumerate}

We fill in the details in the following section, in reverse order:
The calculation in item~\ref{five} is carried out in 
Lemma~\ref{RqpkJk}, to find the sheaves in item~\ref{four}.
Item~\ref{three} is done in Lemma~\ref{HmTkJk}, where we use the 
Leray spectral sequence to find $H^m(T^n;\J')$ by an induction argument.
We apply this to find a description of the sheaf $R^q\pi_*\J$ 
in Corollary~\ref{RqpiJ}.  
In Lemma~\ref{HpVRqpiJ}, we apply the previous result to find 
$H^p(V;R^q\pi_*\J)$, as in~\ref{one}.
And finally, in Theorem~\ref{lerayband} we use this to find 
the sheaf cohomology of a band in $(\R\cross S^1)^n$.

In section~\ref{leray-ndsing}, we carry out the same calculation 
in the case where the leaf is partially singular.  
The idea is the same, but the calculation is simpler, and 
we only use the Leray spectral sequence once, 
applied to the map that projects out one of the singular components.

We remark that this outline shows the purpose of considering the stalks 
of the higher direct image sheaves, rather than trying to compute
the sheaves themselves from Definition~\ref{Rqsheaf}.  
If we were to follow through this calculation considering 
the sheaf over a small open set at each stage instead of the stalk,
we would in step~\ref{five} be computing not 
$H^q(\{\text{pt}\}\cross S^1;\J')$ but $H^q(W\cross S^1;\J)$, where $W$
is some small \nbhd\ in $T^{n-1}\cross\R^n$.
To do this directly by a method similar to section~\ref{s:cylinder},
we would have to use a cover of a $(2n-1)$-dimensional set, 
which becomes combinatorially unmanageable. 
The sheaf theorist's trick permits us to use induction in 
step~\ref{five} instead.

\subsection{The case of several dimensions: non-singular}\labell{leray-nd}

As noted above, let $\J_n$ be the sheaf of leafwise flat sections
over $(\R\cross S^1)^n$.
Let $\J'_n$ denote the restriction of $\J_n$ to $\{\bt\}\cross(S^1)^n$.  
(The sheaf $\J'_n$ will depend on the value of $\bt$, though
this is not made explicit in our notation.)
Elements of $\J'_n$ will have the local form 
\begin{equation}\labell{Jprime}
a'(t_1,\ldots,t_n) e^{i(t_1\theta_1 + \cdots + t_n \theta_n)},
\end{equation}
where $t_1,\ldots t_n$ are constant, and 
 $a'$ is the germ of a smooth function of the $t$ variables.
We will be very careful in the following to 
note what sheaf we are working with.
Let \[ \pi \colon (\R\cross S^1)^n \to \R^n \]
be the projection to the $\R$ factors, and let
\[ \rho \colon T^n \to T^{n-1} \]
be projection to the first $n-1$ factors.

\begin{lemma}  \labell{RqpkJk}
The higher direct images of $\J_n'$ are 
\begin{equation}
R^q \rho_* \J'_n \cong 
\begin{cases}
\J'_{n-1} &\quad \text{if } t_n \in \Z \text{ and } q=1 \\
0 &\quad \text{otherwise}
\end{cases}
\end{equation}
\end{lemma}

\begin{proof}
First, the case $q=0$ is easy to see directly:
Let $U$ be a small open neighbourhood in $T^{n-1}$.  
The preimage $\pi^{-1}(U)$ is 
$U \cross S^1$, and $(\rho_* \J'_n) (U)$ is just the set of 
elements of $\J'_n$ over $U\cross S^1$.  
As noted above in \eqref{Jprime}, elements of $\J'_n$ can be written 
locally as
$a' e^{i(t_1\theta_1 + \cdots + t_{n-1} \theta_{n-1})} e^{it_n\theta_n}$, 
where $(\theta_1, \ldots, \theta_{n-1}) \in U$ and $\theta_n$ ranges 
from 0 to $2\pi$.

In order for such a creature to be defined on the whole set, 
it must have the same value for $\theta_n=0$ as for $\theta_n=2\pi$,
\emph{as germs of functions of the $t$ variables.}
This is impossible, even if $t_n \in \Z$---for example, 
$e^{it\theta}$ and $e^{it\theta + 2\pi it}$ are different as germs, 
even if their values are the same for $t \in \Z$.
Thus there are no elements of $\J'_n$ defined on all of $U\cross S^1$, 
and so $\pi_*\J'_n = 0$.

\vspace{1ex}
Next, consider the case $q=1$.

Recall that 
\begin{equation}\labell{H1thS1}
 (R^1 \rho_*\J'_n)_{\bth} \cong H^1 (\rho^{-1}(\bth),\J'_n) 
	= H^1(\bth \cross S^1, \J'_n) 
\end{equation}
View $\J'_n$ as a sheaf over $S^1$.
This calculation is very similar to the one in 
section~\ref{ss:sheaf-cylinder}; the difference is that we will
be dealing with germs in the $t$ variables, instead of functions of $t$.

Cover $S^1$ with three sets $E$, $F$, and $G$, as in section 
\ref{ss:sheaf-cylinder}.
An element of $\J'_n$ over $S^1$ looks like 
\[ a' e^{i(t_1 \theta_1 + \cdots + t_{n-1} \theta_{n-1})} 
	e^{i t_n \theta_n} 	\]
where $\theta_1, \ldots , \theta_{n-1}$ are fixed, and $\theta_n$ 
ranges over all values from 0 to $2\pi$, 
and all the $t$ variables are fixed.
This element is determined by $a'$, which is the germ of a smooth function
of $(t_1, \ldots, t_n)$.

A \v Cech 1-cochain is a collection of three elements of $\J'_n$, one
for each intersection of $E$, $F$, and $G$.  Thus such a cochain 
is determined by three germs $\{ a'_{EF}, a'_{FG}, a'_{GE} \}$.
All 1-cochains are cocycles.

The calculation of this cohomology is more or less the same as 
in section~\ref{ss:sheaf-cylinder}.  
In order for the cochain $\{ a'_{EF}, a'_{FG}, a'_{GE} \}$ to be a 
coboundary, we need a 1-cochain $\{ b_E, b_F, b_G \}$ which is a primitive.

This condition leads us to a set of equations like~\eqref{1-coc}, 
except with some extra factors:
\begin{subequations}
\begin{align}
a_{EF} e^{it_n\theta_E}e^{i\ip{\bt}{\bth}} 
	&= b_F e^{it_n\theta_F}e^{i\ip{\bt}{\bth}} 
	- b_E e^{it_n\theta_E}e^{i\ip{\bt}{\bth}} \qquad \text{on }E\cap F\\
a_{FG} e^{it_n\theta_F}e^{i\ip{\bt}{\bth}} 
	&= b_G e^{it_n\theta_G}e^{i\ip{\bt}{\bth}} 
	- b_F e^{it_n\theta_F}e^{i\ip{\bt}{\bth}} \qquad \text{on }F\cap G\\
a_{GE} e^{it_n\theta_G}e^{i\ip{\bt}{\bth}}
	 &= b_E e^{it_n\theta_E}e^{i\ip{\bt}{\bth}}
	 - b_G e^{it_n\theta_G}e^{i\ip{\bt}{\bth}} \qquad \text{on }G\cap E
\end{align}
\end{subequations}
Here each of the $a$'s and $b$'s are germs in $t_1,\ldots t_n$.  

The argument following \eqref{1-coc} in section~\ref{ss:sheaf-cylinder}
goes through unchanged, working with germs rather than functions.  
(The extra factor of $e^{i\ip{\bt}{\bth}}$, since it is nonzero, can
be cancelled from each equation.  
It is also worth noting that the role of $t$ in section~\ref{ss:sheaf-cylinder}
is played here by $t_n$.)
We get the same matrix \eqref{matrix} for the system, and the same 
condition that the system has a solution if $e^{-2\pi i t_n} \neq 1$.
We can think of $t_n$ as a parameter of the sheaf, and if it is 
not an integer, then the cohomology in~\eqref{H1thS1} is zero.  

If $e^{-2\pi i t_n} = 1$, then by the same linear algebra argument, 
the system has a solution only if 
\begin{equation}
a_{EF}\restr{t_n} + a_{FG}\restr{t_n} + a_{GE}\restr{t_n} = 0.
\end{equation}
Thus the cohomology \eqref{H1thS1} is given by 
\begin{equation}
H^1(\bth \cross S^1, \J'_n) \cong \{\text{germs in } t_1,\ldots t_n \}
	/ \{ a_{EF}\restr{t_n} + a_{FG}\restr{t_n} + a_{GE}\restr{t_n} = 0 \}.
\end{equation}

Define a map from 1-cochains to $\J'_{n-1}$ by
\[ \{ a'_{EF}, a'_{FG}, a'_{GE} \} \mapsto 
(a'_{EF} + a'_{FG} + a'_{GE}) 
e^{i(t_1 \theta_1 + \cdots + t_{n-1} \theta_{n-1})}. \]
Its kernel is cocycles which have 
$\bigl( a'_{EF} + a'_{FG} + a'_{GE} \bigr) \restr{t_n} = 0$,
i.e.\ coboundaries.  Therefore, 
\[ \J'_{n-1}(U) \cong \{ \text{cocycles} \} / 
\{ \text{coboundaries} \}  \cong H^1(U;\J'_n), \]
which equals $R^1\rho_*\J_n'$.

Finally, a similar technique works for the higher direct images%
---the argument of Section~\ref{ss:sheaf-cylinder},
with germs instead of functions, shows that 
$R^q \rho_* \J'_n = 0$ for $q\geq 2$.
\end{proof}

\begin{lemma}\labell{HmTkJk}
Assume $\bt \in \Z^n$.  Then 
\begin{equation}
H^m(\{\bt\} \cross T^n;\J_n') \cong
\begin{cases}
0 \qquad\text{if } m\neq n\\
\C \qquad\text{if } m=n.
\end{cases}
\end{equation}
If $\bt \notin \Z^n$, then $H^m(\{\bt\} \cross T^n;\J_n') =0$.
\end{lemma}

\begin{proof}
Assume first that all coordinates of $\bt$ are integers.
We proceed by induction on $n$.

First, the case of $n=1$ is just Lemma~\ref{RqpiJ1}: $H^1 \cong \C$, 
$H^\text{other} = 0$.

Consider the Leray spectral sequence for the map $\rho$.
It will have $E_2$ term given by 
\[ E_2^{p,q} = H^p(T^{n-1};R^q \rho_* \J'_n). \]

By Lemma~\ref{RqpkJk}, all of the $R^q\rho_*\J'_n$ are zero
except $R^1$.  
Therefore, only one row of the spectral sequence is non-zero, 
and thus, as discussed in Proposition~\ref{ss-triv}, the 
spectral sequence immediately stabilizes and $E_\infty = E_2$.
Thus 
\begin{equation}\labell{sumS1Rqpk}
\begin{split}
 H^m(T^n;\J'_n) = E_\infty^m 
	&\cong \bigoplus_{p+q=m} H^p(T^{n-1};R^q\rho_*\J_n)\\
	&\cong H^{m-1}(T^{n-1};R^1\rho_*\J_n')
\end{split}
\end{equation}
where the last line is true because all $R^q$ are zero except $R^1$.
Also by Lemma~\ref{RqpkJk}, $R^1\rho_*\J'_n \cong \J'_{n-1}$.  
Thus 
\[ H^m(T^n;\J'_n) \cong H^{m-1}(T^{n-1};\J'_{n-1}). \]
The induction is complete.

\vspace{1ex}

Now, suppose that $\bt \notin \Z^n$.
Without loss of generality, suppose $t_n \notin \Z$.  
Then by Lemma~\ref{RqpkJk}, \emph{all} the sheaves $R^q\rho_*\J'_n$ 
are zero, and so all summands in \eqref{sumS1Rqpk} are zero.
\end{proof}

\begin{cor}\labell{RqpiJ}
The sheaf $R^q \pi_* \J_n$ over $\R^n$ is 0 for $q\neq n$.  
When $q=n$, it is a skyscraper sheaf
supported on $\Z^n$, with each tower isomorphic to $\C$.
\end{cor}

\begin{proof}
Assume first that $\bt \in \Z^n$.

By the sheaf theorist's trick, a stalk of $R^q \pi_* \J_n$ has the form
\[ (R^q \pi_* \J_n)_\bt \cong  
	H^q \bigl( \pi^{-1}(\bt ), \J'_n \bigr)
	\cong H^q \bigl( T^n, \J'_n \bigr). \]
In Lemma~\ref{HmTkJk}, we just showed that this 
is isomorphic to $\C$, if $q=n$, and 0 otherwise.
Therefore, the sheaf $R^1 \pi_* \J_n$ has stalk $\C$ when $\bt \in \Z^n$, 
and zero otherwise.
\end{proof}

\begin{lemma}\labell{HpVRqpiJ}
For $V \subset \R^n$, 
\begin{equation}
H^p(V,R^q\pi_*\J) = 
\begin{cases}
0 \quad \text{if } q \neq n \\
\C^\nu \quad \text{if } p=0 \text { and } q=n \\
0 \quad \text{if } p \geq 1
\end{cases}
\end{equation}
Here $\nu$ is the number of points in the intersection $V \cap \Z^n$.
\end{lemma}

\begin{proof}
This all follows from Corollary~\ref{RqpiJ}.

The first case follows since $R^q \pi_* \J = 0$ if $q \neq n$, 
and so the cohomology is zero 
(the cohomology of a zero sheaf is zero).

The second case follows since the sheaf is a skyscraper.  
As given in Lemma~\ref{l:skys},
the $0^\text{th}$ cohomology of a skyscraper sheaf is the direct 
sum of the towers.
Since $R^n \pi_* \J$ has tower $\C$ over each point of $\Z^n$, 
the result for $p=0$, $q=n$ follows.

The third case follows because the higher cohomology of a skyscraper
sheaf is zero.
\end{proof}

\begin{thm}\labell{lerayband}
For $U\subset (\R\cross S^1)^n$ a band,
$H^n(U,\J) \cong \C^\nu$, where $\nu$ is the number of
Bohr-Sommerfeld orbits contained in $U$.
For all other values of $m$, $H^m(U,\J) = 0$.
\end{thm}

\begin{proof}
Let $\nu$ be the number of Bohr-Sommerfeld leaves contained in $U$, 
and let $V = \pi(U) \subset \R^n$.  
Since the \BS\ leaves are precisely the fibres of $\pi$ over 
the points in $\Z^n$, $\nu$ is the number of points 
in the intersection $V \cap \Z^n$.

At last, we will apply the Leray spectral sequence to the map 
$\pi \colon (\R\cross S^1)^n \to \R^n.$
The $E_2$ term of the sequence is 
\begin{equation}
\E_2^{p,q} = H^p(V;R^q\pi_*\J). 
\end{equation}
According to Lemma~\ref{HpVRqpiJ}, this is only non-zero if 
$(p,q) = (0,n)$.
Thus the spectral sequence has only one non-zero entry 
(never mind one non-zero row), 
and so by Proposition~\ref{ss-triv}, the spectral sequence
stabilizes immediately and 
\begin{equation}
H^m(U,\J) = \bigoplus_{p+q=m}  H^p(V,R^q\pi_*\J). 
\end{equation}
The only non-zero summand is the $(0,n)$ one, which again by 
Lemma~\ref{HpVRqpiJ} is isomorphic to $\C^\nu$.
\end{proof}


\subsection{The partially singular case}\labell{leray-ndsing}

We now deal with the case where the Bohr-Sommerfeld leaf is singular.
This means that the model space is $(\R\cross S^1)^{n-k} \cross \C^k$, 
for some $k>0$.  
In this case, a simple neighbourhood of the leaf is a more general
band: 
the product of bands around the cylinder components, and a product 
of discs in the $\C$ components (see Definition~\ref{def:band}).

\begin{thm}
If $U \subset (\R\cross S^1)^{n-k} \cross \C^k$ 
is a band, and $k>0$, 
then $H^m(U;\J) = 0$ for all values of $m$.
\end{thm}

\begin{proof}
Assume, for definiteness, that $k=1$, so 
$U \subset (\R\cross S^1)^{n-1} \cross \C$. 

Let 
\[p \colon U \to Y=(\R\cross S^1)^{n-1} \] 
be projection onto the first $n-1$ coordinates.
The Leray spectral sequence for the map $p$ gives 
$H^m(U;\J)$ in terms of $H^r(Y;R^qp_*\J_n)$.
By the sheaf theorist's trick (Theorem~\ref{stt}), 
the sheaf $H^r(Y;R^qp_*\J_n)$ is given by
\begin{equation}
(R^qp_*\J)_x \cong H^q(p^{-1}(x);\J_n'),
\end{equation}
where $\J_n'$ is the sheaf $\J_n$ restricted to 
$p^{-1}(x) = \{x\} \cross D^2$.

Elements of the sheaf $\J$, recall, have the local form 
$a(t,s) e^{i \ip{t}{\theta}} e^{i\ip{s}{\phi}}$.
Since we are assuming for the moment that $k=1$, 
elements of $\J'$
locally have the form $a'(\bt,s) e^{i\ip{\bt}{\bth}} e^{is\phi}$, 
where $a'$ is the germ of a smooth function of all the $(n-1)$ $t$ variables, 
and a smooth function of $s$.  
In $D^2$, the $\phi$ and $s$ variables can change; 
the $\bt$ and $\bth$ variables are fixed by $x$ 
(in fact $x = (\bt,\bth)$).
Call such a creature a ``semigerm''.  
 Thus we wish to calculate 
\[ H^q ( \{x\} \cross D^2; \J_n' ) \]
where $\J_n'$ is the sheaf of semigerms.

This calculation is almost identical to the one in Section~\ref{C-section},
except that the elements of the sheaf are of the form 
$a'(\bt,s) e^{i\ip{\bt}{\bth}} e^{is\phi}$
instead of $a(s) e^{is\phi}$.  
The same calculation as the one following Proposition~\ref{prop-annulus}
goes through in this case, except that the coefficients 
$h_{EF}$, etc., are semigerms instead of just functions of $s$,
and there is an extra factor of $e^{i\ip{\bt}{\bth}}$ multiplying everything.
For example, equation \eqref{pre-cocyc} becomes
\begin{equation}
\begin{split}
h_{EF} e^{\ip{\bt}{\bth}} e^{is\phi_E} - h_{AF} e^{\ip{\bt}{\bth}} e^{is\phi_F}
	+ h_{AE} e^{\ip{\bt}{\bth}} e^{is\phi_E} &= 0 \qquad \text{on } AEF\\
h_{FG} e^{\ip{\bt}{\bth}} e^{is\phi_F} - h_{AG} e^{\ip{\bt}{\bth}} e^{is\phi_G}
	+ h_{AF} e^{\ip{\bt}{\bth}} e^{is\phi_F} &= 0 \qquad \text{on } AFG\\
h_{GE} e^{\ip{\bt}{\bth}} e^{is\phi_G} - h_{AE} e^{\ip{\bt}{\bth}} e^{is\phi_E}
	+ h_{AG} e^{\ip{\bt}{\bth}} e^{is\phi_G} &= 0 \qquad \text{on } AGE
\end{split}
\end{equation}
where each $h_{UV} (\bt,s)$ is 
a germ in the $t$ variables and a smooth function of $s$.

Regardless of the (fixed) values of $\bt$ and $\bth$, 
$ e^{\ip{\bt}{\bth}}$ will be nonzero, and so 
it can be cancelled from each of the equations.  
Also, as in section~\ref{ss:C-cohom}, we can cancel the factors
of $e^{is\phi}$, keeping track of the conventions on the $\phi$ coordinates
(given in~\ref{ss:C-cohom} just before equation \eqref{cocycle}).
This gives a set of equations identical to \eqref{cocycle},
this time for the semigerms $h_{EF}$, $h_{FG}$, etc., and
the calculation proceeds exactly the same, 
and gives the same result for $\{x\}\cross D^2$ as in 
Theorem~\ref{disczero}, that the degree 1 cohomology is zero.

Similarly, the argument in the proof of Proposition~\ref{p:C-cohom-k} 
goes through unchanged, with semigerms and the 
extra factor of $e^{\ip{\bt}{\bth}}$,
and tells us that the cohomology in all other degrees is trivial.

Finally, we return to the calculation of $H^m(U;\J)$.
By what we have just shown, $R^q p_* \J$ is zero for all $q\geq 0$.
Thus $H^r(Y;R^qp_*\J)=0$ for any choice of $r$ and $q$, and so 
the Leray spectral sequence for the map $p$ not only stabilizes, 
it has all $E_2$ terms equal to zero.  Therefore
$H^m(U;\J) =0$ for all $m$.

Note we have assumed for simplicity that $k=1$, but the same argument will
hold if $k>1$ as well; there will just be more $s$ and $\phi$ coordinates.
Thus 
\begin{equation}
H^m(U;\J) =0 \qquad \text{for all } m
\end{equation}
whenever the leaf is at all singular.
\end{proof}

\section{Piecing and glueing}\labell{s:piecing}

The conditions in Definitions~\ref{spacedef} guarantee that our spaces of 
interest have neighbourhoods which look like open sets in the model spaces.
We will use the results on the cohomology of the model spaces to 
obtain results about the cohomology of these spaces.
In this section we develop the theory necessary to transfer results 
about sheaf cohomology from one manifold to another.

The fact that our spaces are locally symplectomorphic to $\msp$ 
is a part of the definition.  
In section~\ref{ss:fsharp} we also obtain, 
with no extra hypotheses, a stronger condition on the symplectomorphism, 
which will enable us to compare sheaf cohomology.
Finally, in \ref{s:patch}, we put everything together to prove 
our main theorem about the cohomology of thse spaces.

First, though, we need some theory.

\subsection{Necessary sheaf theory}

\begin{lemma}\labell{l:fsharp}
Let $\sh$ and $\T$ be sheaves over manifolds $M$ and $N$, respectively,
and $f \colon M \to N$ a smooth map.  
Suppose that, for every open set $V\subseteq N$, we have a group 
homomorphism\footnote{In our applications, $f^\sharp$ will be 
induced by $f$, hence the notation.}
\[ f^\sharp \colon \T(V) \to \sh\bigl( f^{-1}(V) \bigr) \]
which is compatible with both $f$ and the restriction maps in the sheaves,
in the sense that
\begin{equation}\labell{eq:fsharp}
(f^\sharp \alpha)\restr{f^{-1}(V)} = f^\sharp(\alpha\restr{V})
\qquad \text{for all } V \text{ and } \alpha \in \T(V).
\end{equation}
Then, for any open $V \subseteq N$,  $f^\sharp$ induces a map 
\[ f^* \colon  H^*\bigl(f^{-1}(V),\T\bigr) \to H^*(V,\sh). \]
Furthermore, this process is `functorial,' in the sense that 
the composition of induced maps is the map induced by the 
composition.\footnote{ 
In more detail: Suppose we have maps 
\[ M \stackrel{f}{\to} N \stackrel {g}{\to} P, \]
and maps on sheaves 
\[ \Q \stackrel{g^\sharp}{\to} \T \stackrel{f^\sharp}{\to} \sh. \]
Suppose further that $h=g\circ f$, $h^\sharp = f^\sharp \circ g^\sharp$,
and all the requisite compatibilities are satisfied.  
Then $h^* = f^* g^*$ as maps on sheaf cohomology.}
\end{lemma}

\begin{proof}
In brief, the map $f^\sharp$ induces a cochain map which is compatible with 
the coboundary, and thus induces a map on cohomology.  
The details are straightforward, and are left to the reader.
%
\end{proof}

\subsection{The induced map on cohomology}\labell{ss:fsharp}

\begin{thm}\labell{fsharp}
Let $N$ be a compact symplectic
manifold with prequantization line bundle $(\LL_N, \nabla^N)$,  
equipped with a locally toric singular Lagrangian fibration.
Let $V$ be a \nbhd\ of a leaf $\lf_N$
symplectomorphic to a \nbhd\ $U$ of a leaf $\lf_0$ in a model space 
$M_0 = (\R\cross S^1)^{n-k} \cross \C^k$
(the existence of $V$ and $U$ is guaranteed by Definition~\ref{spacedef}).
Then there exists an invertible map 
$f^\sharp \colon \J_{M_0}\restr{U} \to \J_N\restr{V}$.
\end{thm}

The proof of Theorem~\ref{fsharp} proceeds by a series of lemmas.
Definition~\ref{spacedef} guarantees the existence of a symplectomorphism
$f\colon U \to V$.
We first find a \triv\ of $\LL_0$, 
with respect to which  its connection has 
potential one-form $f^* \Theta_N$ 
(where $\Theta_N$ is the potential one-form of $\nabla^N$).
Then we use this \triv\ to define $f^\sharp$ in such a way
that it takes flat sections to flat sections,
proved in Lemma~\ref{pbflat}.

\begin{lemma}\labell{l:hol}
We may choose $\lf_0$ so that the holonomies of $\lf_0$ and $\lf_N$ 
are equal.
\end{lemma}

\begin{proof}
As noted above, we already have a symplectomorphism $f\colon U \to V$.
Changing $t$ (the coordinate on $\R^{n-k}$) 
by a constant doesn't change the symplectic form on $M_0$,
and so we are free to choose $t_0$, the value of $t$ 
corresponding to $\lf_0= f^{-1}(\lf_N)$.

Each leaf of the singular fibration is homeomorphic to $T^{n-k}$
for some $k$---this is another consequence of the definition.  
Write $m=n-k$, as before.
Let $\beta_1,\ldots,\beta_m $ be the set of fundamental cycles in $T^m$ 
given by loops around each $\theta$ coordinate.  
These are then mapped by $f$ to a set of fundamental cycles for $\lf_N$, 
which we denote $\gamma_1,\ldots,\gamma_m$.

The holonomy around $\gamma_j$ is 
\[ \exp \left( i \int_{\gamma_j} \Theta_N \right) \]
(see section~\ref{s:hol}),
while the holonomy around $\beta_j$ is 
\[ \exp \left( i \int_{\beta_j} \Theta_0 \right) = 
 \exp \left( i \int_{\beta_j} 
 \sum
 \bigl( t_a \, d\theta_a + s_b \, d\phi_b \bigr)
 \right) = e^{2\pi it_j}.
\]
Thus, if we choose the $j^\text{th}$ coordinate of $t_0$ to be 
\[ \frac{1}{2\pi} \int_{\gamma_j} \Theta_N \]
for each $j=1,\ldots,m$, then the holonomies of $\lf_N$ and $\lf_0$ 
will be equal.
\end{proof}

\begin{lemma}
The bundle $\LL_N$ is trivializable over $V \subset N$.
\end{lemma}

\begin{proof}
The hypotheses on our spaces guarantee that 
we can choose 
$V$ to be of the form $I^m \cross T^m \cross (D^2)^k$,
with the leaf $\lf_N$ being identified with the central torus
$\lf_N \cong \{ \bt^m \} \cross T^m \cross \{ 0 \}$.
A transverse \nbhd\ is just $I^m \cross D^{2k}$, which is contractible; 
therefore $\LL_N$ is trivializable over it.  
There is a free $T^m$ action on $V$ (just act on the $T^m$ coordinate),
which ``sweeps out'' the transverse disc over the \nbhd\ $V$.
This action gives us a \triv\ over the whole \nbhd. 
\end{proof}

\begin{lemma}\labell{step2}
Choose a local \triv\ of $\LL_N$ over $V$, and say 
$\Theta_N$ is the potential one-form of $\nabla^N$ with respect 
to this trivialization.  
Then there is a \triv\ of $\LL_0$ over $U$ with respect
to which $\nabla^0$ has potential one-form $f^* \Theta_N$.
\end{lemma}

\begin{proof}

Let $\Theta_1 = f\pb \Theta_N$ be the pullback of the potential one-form.
We wish to show that there is a trivialization of $\LL_0$ with respect 
to which the connection on $M_0$ has potential one-form $\Theta_1$.

Since the curvature of each connection is the symplectic 
form on the respective manifold, and $f$ is a symplectomorphism, 
$d\Theta_0 = f\pb d\Theta_N = d\Theta_1$.  
Therefore $\Theta_1 - \Theta_0$ is closed.

Recall from the proof of Lemma~\ref{l:hol} that we have 
chosen $t_0$ so that 
\[ \int_{\gamma_j} \Theta_N = \int_{\beta_j} \Theta_0 \]
for all $j$.  
Pushing the left side forward by $f^{-1}$, we obtain that 
\[ \int_{\beta_j} \Theta_1 = \int_{\beta_j} \Theta_0 \]
for all $j$.
Therefore $\Theta_1 - \Theta_0$ is closed and has integral 0 around 
all loops in $T^k$, and thus in $\lf_0$, and thus in $U$.
Therefore, it is exact.
(Note that the $(s,\phi)$ coordinates do not enter into this 
consideration, since $U$ is a disc in those coordinates, 
and so there are no nontrivial loops in the $(s,\phi)$ coordinates.)
Write $\Theta_1 = \Theta_0 + dG$
with $G\colon U \to \R$.

Now, suppose we change the trivialization of $\LL_0$ by multiplying 
the fibres by some (nonzero) function $\psi\colon U \to S^1$; 
i.e., we take $(p, z) \mapsto (p, \psi(p) z )$.
By Eq \eqref{eq:trivchg}, this changes the potential 
one-form of a connection by subtracting $i \tfrac{1}{\psi}\, d\psi$ from it.
Thus, we require  a function $\psi$ 
such that 
\[ i \tfrac{1}{\psi}\, d\psi = - dG. \]
Such a function is 
\[ \psi = e^{i G}. \]
Thus, multiplying the canonical \triv\ of $\LL_0$ by $e^{iG}$ gives 
a new \triv, with respect to which the connection
$\nabla^0$ has potential one-form $\Theta_1$.
\end{proof}

\begin{defn}
Given $N$, $M_0$, etc.\ as above, 
define the map $f^\sharp \colon \Gamma(V, \LL_N) \to \Gamma(U, \LL_M)$
as follows.  Let $s_V$ be the unit section of $\LL_N$ in some 
trivialization over $V$, and $\Theta_N$ the corresponding 
potential one-form of $\nabla^N$.
By Lemma~\ref{step2} there is a trivializing section $s_U$ of $\LL_0$ 
over $U$, with respect
to which $\nabla^0$ has potential one-form $f\pb \Theta_N$.
If $\sigma$ is a section of $\LL_N$, then $\sigma = \phi\, s_V$ for a 
function $\phi$.
Then $f^\sharp \sigma$ is the section $(\phi \circ f)\, s_U$.
\end{defn}

Note that $f^\sharp$ is invertible, since $f$ is invertible:
map $\phi\, s_U$ to $(\phi \circ f^{-1})\, s_V$.

\begin{lemma} \labell{pbflat}
If $\sigma$ is flat along the leaves of $N$, then 
$f^\sharp \sigma$ is flat along the leaves of $M_0$.
\end{lemma}

\begin{proof}
Since $f$ is a diffeomorphism which carries leaves to leaves, it will 
suffice to prove that $\nabla^M_X f^\sharp \sigma =0$ whenever 
$\nabla^N_{f_*X} \sigma = 0$.  
This follows from chasing the definitions, using a couple 
of facts about pullbacks and the chain rule.
\end{proof}

\begin{cor}
The map $f^\sharp$ defined above is a sheaf map from $\J_{M_0}\restr{U}$
to $\J_N\restr{V}$.
\end{cor}

\begin{proof}
Clearly $f^\sharp$ is compatible with the restriction maps, 
and thus it is a sheaf map.  Since it takes flat sections to 
flat sections, it maps $\J_{M_0}$ to $\J_N$.
\end{proof}

Given $N$, we have constructed an invertible map from
$\J_{M_0}$ to $\J_N$.  Thus the  proof of Theorem~\ref{fsharp} is complete.

\begin{cor}\labell{sheafiso}
Under the conditions of Theorem~\ref{fsharp}, 
\[ H^* (U,\J_M) \cong H^*(V,\J_N). \]
\end{cor}

\begin{proof}
By Theorem~\ref{fsharp}, $f^\sharp$ satisfies the conditions of 
Lemma~\ref{l:fsharp}, and thus induces a homomorphism $f^*$ on cohomology.
Note that $f^\sharp$ is invertible, and its inverse $g^\sharp$ 
satisfies the functoriality conditions in the footnote to Lemma~\ref{l:fsharp}.
Therefore, $f^*$ is an isomorphism, by the following standard argument:
The maps induced on cohomology by $f^\sharp$ and $g^\sharp$ satisfy
\[ f^* g^* = \id; \qquad g^* f^* = \id \]
by functoriality, and so they must be isomorphisms.
\end{proof}

\subsection{Patching together}\labell{s:patch}

Finally, we are in a position to state and prove our main theorem.  
First a small lemma.

\begin{lemma}\labell{pieces}
 Let $M$ be a compact $2n$-dimensional symplectic manifold, 
equipped with a locally toric singular Lagrangian fibration.
Then $M$ can be covered by finitely many open bands,
such that any \BS\ leaf is contained in only one band.
\end{lemma}

\begin{proof}
Recall first the 
definition of our spaces of interest, Definition~\ref{spacedef}.

Cover the base $B$ by sets $U$ which are homeomorphic to 
open rectangles 
in $\R^{n-k} \cross \R_+^k$.  By shrinking them, if necessary, 
we may assume that no \BS\ point lies in more than one of them.
Then the inverse images of the $U$'s in $M$ will be symplectomorphic
to the inverse images of the rectanges in $\R^{n-k} \cross \R_+^k$,
which is precisely the definition of a band in this context 
(Definition~\ref{def:band}), and no \BS\ leaf will lie in 
more than one of them.
Finally, finitely many of them will suffice to cover $M$ since $M$ is compact.
\end{proof}

\begin{thm}[Main Theorem]\labell{mainthm}
Let $M$ be a compact $2n$-dimensional symplectic manifold, 
with a prequantization line bundle $\LL$, and 
with a (singular) real polarization
given by a locally toric singular Lagrangian fibration.
Let $\J$ be the sheaf of leafwise flat sections of $\LL$.
Then the cohomology groups $H^q(M;\J)$ are zero for all $q \neq n$, and 
\begin{equation}\labell{M-cohom}
 H^n (M;\J) \cong \bigoplus_{b \in BS} \C 
\end{equation}
where the sum is taken over all non-singular \BS\ fibres.
\end{thm}

\begin{proof}
First, cover $M$ by sets $U$ as in Lemma~\ref{pieces}.
As each $U$ is equivalent to a band in $\msp$, 
via a symplectomorphism satisfying the conditions of Theorem~\ref{fsharp},
it has the same cohomology as such a generalized band.  
By the results in section~\ref{s:leray}, a band has sheaf cohomology 
$H^q = 0$ for all $q\neq n$, and its $n^\text{th}$ 
cohomology has one copy of $\C$ for each non-singular \BS\ leaf 
it contains.
By Mayer-Vietoris (Prop~\ref{mult-MV}), the cohomology of $M$ is 
the sum of the cohomology of each of these generalized bands in the cover
(since we are assuming no \BS\ leaf lies in the intersection of any two 
sets in the cover),
which gives us \eqref{M-cohom}.
\end{proof}

\section{Real and K\"ahler polarizations compared}\labell{s:toric}

As noted in section~\ref{ss:qn-pol}, the (singular) foliation of 
a toric manifold by the fibres of the moment map is a singular 
Lagrangian fibration, and thus is a singular real \poln .
Thus, by the above theorem, the \qn\ of a toric manifold has dimension 
equal to the number of nonsingular \BS\ leaves.  

\begin{prop}\labell{intlattice}
For $M$ a toric manifold with moment map 
$\mu \colon M \to \Delta \subset \R^n$, 
the Bohr-Sommerfeld set is the set of integer lattice points in the 
moment polytope $\Delta$.
The singular \BS\ points are the ones on the boundary of the polytope.
\end{prop}

(This is a well-known result, but we include a proof here for the 
sake of completeness.)

\begin{proof}
Guillemin and Sternberg discuss the connection between Bohr-Sommerfeld 
points and action-angle variables in section 2 of~\cite{GS}.
Given a Lagrangian fibration $\pi\colon M \to B$ with compact fibres,
they construct action coordinates as follows: 
Assume that $p \in B$ is a regular value of $\pi$.
Choose a neighbourhood 
$V\subset B$ of $p$ such that $\omega$ is exact on $\pi\inv(V)$, 
with $\beta$ a primitive for $\omega$.  
The fibres of $\pi$ for points in $V$ are tori; 
choose a homology basis 
$\gamma_1(q), \ldots, \gamma_n(q)$ for the fibre over $q\in V$ which 
depends continuously on $q$.  
Define the action coordinates $a_j \colon V \to \R$ by 
\[ a_j(q) = \frac{1}{2\pi} \int_{\gamma_j(q)} \beta. \]
The $\{ a_j \}$ are defined up to addition of a constant 
(which comes from changing $\beta$) 
and multiplication by an element of $GL(n,\Z)$
(from changing the homology basis $\{ \gamma_j \}$).

The holonomy around $\gamma_j(p)$ and $\gamma_j(q)$ 
differ by $\exp i\bigl( a_j(p) - a_j(q) \bigr)$.  
Therefore, if the action coordinates are normalized so that 
$a_j(p) = 0$ at some \BS\ point $p$, then another point 
$q$ is a \BS\ point if{f} all $a_j(q)$ are integers.

For a toric manifold, the coordinate system consisting of coordinates 
on $\Delta \subset \R^n$ together with coordinates on the torus fibre
give action-angle coordinates, more or less by definition.  
Furthermore, because $\Delta$ is a polytope and thus simply-connected,
these action-angle coordinates can be defined on all of
the interior of the polytope.  
The discussion in~\cite{GS} is for the case of a fibration, which here only 
applies to the regular values of $\mu$, 
but it is not hard to extend the result to the singular values, as follows.

Assume that $p_0$ is a fixed point of the action; the case for 
a more general singular point is similar.
By the Local Normal Form, we have a system of coordinates 
$(s_j, \phi_j)$ on a nieghbourhood of the origin in $\C^n$,
as in the local model space 
in section~\ref{s:multdm}.  Points with all $s_j$ nonzero 
are regular points of $\mu$, and
the curves $\gamma_j$ consisting of circles with fixed $s_j$ 
in the $(s_j,\phi_j)$ plane
form a homology basis for the regular fibres.
Thus the $s_j$ are action coordinates on the regular fibres, 
which extend continuously to the fixed point.  

Assume that there is a \BS\ point $p$ in the neighbourhood on which these 
$(s,\phi)$ coordinates are defined.  
By the same argument as in section~\ref{C-section}, all the $s_j(p)$ 
must be integers.  The action coordinates $t_j$ are also defined at 
$p$, since they are defined on the entire interior of $\Delta$, and 
have integer values at $p$.
Therefore (as in (Eq 2.5) in~\cite{GS}),
\[ t_j = \sum a_{ij} (s_i + c_i) \]
for some $c\in \Z^n$ and $A=(a_{ij}) \in GL(n,\Z)$.
The above formula is valid when $s_j \neq 0$ for all $j$, but 
extends continuously to where $s_j=0$, and so extends the definition 
of the action coordinates to all of $M$.
Also, since $c\in \Z^n$ and $A=(a_{ij}) \in GL(n,\Z)$, 
all the $s_j$ are integers if{f} all the $t_j$ are integers.
Therefore the \BS\ points in $\Delta$ are the integer lattice points.

Note that this gives another illustration of the fact 
(mentioned for example in Example 6.10 in~\cite{ggk}) 
that if $M$ is quantizable, 
the fixed points are mapped to integer lattice points in $\R^n$.

Finally, as noted in section~\ref{ss:toric}, 
if $x\in \Delta$ lies on a face of codimension $k$, 
then $\mu\inv(x)$ is an orbit of dimension $n-k$; 
thus, singular orbits correspond to points on the boundary of the polytope.
\end{proof}

A toric manifold also has a natural complex structure, 
coming from its construction as a toric variety, 
and thus a natural K\"ahler \poln .
If the manifold is quantized using this \poln , 
the dimension of the quantization is equal to the number of 
integer lattice points in the moment polytope, 
\emph{including} points on the boundary. 
(This is a well-known result; see~\cite{mdh} for a more complete 
discussion.)

Thus the \qn s coming from K\"ahler and real \poln s
are different, the difference being the number of 
lattice points lying on the boundary of the moment polytope
(which is always non-zero, since fixed points map to lattice points).

\end{document}